\documentclass[%
 aip,
 amsmath,amssymb,
 reprint,%
]{revtex4-1}

\usepackage{graphicx}
\usepackage{dcolumn}
\usepackage{bm}

\usepackage{float}
\usepackage{subfigure}
\usepackage[utf8]{inputenc}
\usepackage{mathptmx}
\usepackage{etoolbox}
\usepackage[export]{adjustbox}
\usepackage{float}
\usepackage{color}
\usepackage{algpseudocode}
\usepackage{algorithm}

\makeatletter
\def\@email#1#2{%
 \endgroup
 \patchcmd{\titleblock@produce}
  {\frontmatter@RRAPformat}
  {\frontmatter@RRAPformat{\produce@RRAP{*#1\href{mailto:#2}{#2}}}\frontmatter@RRAPformat}
  {}{}
}%
\makeatother

\begin{document}

\preprint{AIP/123-QED}

\title[]{Online data-driven changepoint detection for high-dimensional dynamical systems}
\author{Sen Lin}
\altaffiliation[Work performed as ]{Givens Associate at Argonne National Laboratory.}
\affiliation{ 
Department of Mathematics, University of Houston
}
\author{Gianmarco Mengaldo}
\affiliation{Department of Mechanical Engineering, National University of Singapore}
\author{Romit Maulik}%
 \email{rmaulik@anl.gov}
\affiliation{ 
Assistant Computational Scientist, Mathematics and Computer Science Division, 
Argonne National Laboratory,
Lemont, Illinois, 60439, USA.
}%


\date{\today}

\begin{abstract}
The detection of anomalies or transitions in complex dynamical systems is of critical importance to various applications. In this study, we propose the use of machine learning to detect changepoints for high-dimensional dynamical systems. Here, changepoints indicate instances in time when the underlying dynamical system has a fundamentally different characteristic - which may be due to a change in the model parameters or due to intermittent phenomena arising from the same model. We propose two complementary approaches to achieve this, with the first devised using arguments from probabilistic unsupervised learning and the latter devised using supervised deep learning. Our emphasis is also on detection for high-dimensional dynamical systems, for which we introduce the use of dimensionality reduction techniques to accelerate the deployment of transition detection algorithms. Our experiments demonstrate that transitions can be detected efficiently, in real-time, for the two-dimensional forced Kolmogorov flow, which is characterized by anomalous regimes in phase space where dynamics are perturbed off the attractor at uneven intervals.
\end{abstract}

\maketitle

\section{Introduction and motivation}

Changepoint detection is important for dynamical systems because it can help identify shifts or changes in the behavior of a system over time. In many real-world scenarios, it is common for dynamical systems to exhibit non-stationary behavior, i.e., the statistical properties of the system can change over time. Examples of fields where dynamical systems exhibit non-stationary behavior include finance, weather and climate, and fluid dynamics. In the context of fluid dynamics, changepoint detection can be useful for pinpointing regime changes in the dynamics of a fluid flow, involving e.g., velocity, pressure, and/or temperature of the fluid, along with their interactions. Detecting these changes can be important for a range of applications, such as understanding whether the flow is on the verge of instability and turbulence.. Changepoint detection is also crucial in a variety of sectors such as in medical health monitoring for preventing adverse patient outcomes, speech recognition, and human activity detection. \textcolor{red}{\cite{Stival2022,Liu2018,Chowdhury2012, Jung2021}}

Online changepoint detection is a technique that can be used to detect changepoints in real-time data streams. This is in contrast with offline approaches, where a stream of data is post-processed to detect for anomalies after data-collection. Online changepoint detection relies on either a model update, or a model evaluation with new data as it becomes available. Therefore, this approach can adapt to changes in the underlying dynamics of the system, ideally without the use of labelled data. There are several approaches for changepoint detection in dynamical data, each with its own strengths and weaknesses. They may be loosely grouped under

\begin{itemize}
    \item Statistical approaches \cite{fearnhead07,fearnhead19,BOCPD}: Statistical methods are often used to detect changepoints by modeling the data as a sequence of independent random variables. One common approach is to use hypothesis testing to compare the distributions of the data before and after a potential changepoint. Another approach is to use sequential Bayesian methods to update the probability of a changepoint as new data is observed. 
    \item Model-based approaches \cite{Kalouptsidis2011,Srinivasan2019,Chen2023}: Model-based methods assume that the data is generated by a specific dynamical system, and changepoints are detected when the parameters of the model change. For example, one could use a linear dynamical system model and detect changepoints when the system transitions from one set of parameters to another. 
    \item Machine learning approaches\cite{Kotsiantis2007, Munir2018, Pang2020,Ahad2021}: Machine learning methods can also be used for changepoint detection in dynamical data. These approaches typically involve training a classifier to distinguish between different regimes of the dynamical system, and using the classifier to identify the times when the system transitions from one regime to another.
\end{itemize}

In practical applications, the choice of approach depends on the specific characteristics of the dynamical data, such as the complexity of the system, the amount of available data, and the desired level of sensitivity and specificity for detecting changepoints, in addition to the computational budget available during the online phase of detection.

In this work, we focus on online changepoint detection, and we propose:
\begin{itemize}
    \item A dimensionality reduction strategy for high-dimensional dynamical systems that can preserve independence in the reduced-space by using beta-variational autoencoders; 
    \item An unsupervised transition detection algorithm is deployed by using the Bayesian changepoint detection technique, which models the probability of a changepoint in the dynamics in an online manner;
    \item A supervised learning approach based on the use of long short-term memory neural networks for predicting transition events in a streaming manner.
\end{itemize}
The novel methodological developments are deployed for detecting bursting events in the two-dimensional forced Kolmogorov flow. 

The rest of the paper is organized as follows. In section~\ref{sec:related}, we present related and competing work. In section~\ref{sec:dataset}, we introduce the dataset used for testing the novel methods, that are introduced in section~\ref{sec:methods}. In section~\ref{sec:hyperparam}, we introduce the hyperparameter tuning performed for the tests. In section~\ref{sec:results}, we outline the results. In section~\ref{sec:conclusions}, we draw the key concluding remarks of this work.

\section{Related work}\label{sec:related}

Various approaches for changepoint detection have been proposed for different applications in previous literature \cite{aminikhanghahi17,vandenBurg2022}. Let $\boldsymbol{z}_t\in\mathbb{R}$ be independent observations at time $t$, where $t=1,2,...,T$. Let there also exist a set of changepoints $\tau = \{\tau_1,...,\tau_n\}$. Let us denote a segment of the series from $t=a,a+1,...,b$ as $z_{a:b}$. Offline algorithms \cite{Scott1974,Auger1989,jackson05,killick12,fearnhead19} for multiple changepoint detection are developed by optimizing the loss function:
\begin{align}
    \min_{\tau}\sum_{i=1}^n l(z_{\tau_{i}:\tau_{i+1}-1})+\lambda P(n)
\end{align}
Where $l$ is an additive loss function such that $l(z_{a:b})=l(z_{a:\tau-1})+l(z_{\tau:b})$, $P(n)$ is a penalty function on the number of changepoints, and $\lambda\leq0$ is a hyperparameter. In contrast,  online algorithms\cite{fearnhead07,BOCPD} rely on Bayesian posterior probability calculations based on whether the new arrival data is a changepoint. Knoblauch and Damoulas \cite{knoblauch18} added support for model selection and spatio-temporal models, as well as robust detection using $\beta$-divergences. According to Burg and Williams\cite{vandenBurg2022}, offline detection methods are generally restricted to univariate time series while the online method support multidimensional data.

In contrast, model-based approaches rely on online model identification techniques, which sequentially determines model structure and estimates/updates model parameters when new observations arrive\cite{Rong2006,Lombaerts2009,Kopsinis2010,Kalouptsidis2011,Chen2014,Srinivasan2019}. A recent exemplar of this approach is seen in Chen et al., \cite{Chen2023} who proposed the use of a causation entropy boosting (CEBoosting) strategy, which is incorporated into an online model identification method to detect regime switching.

Finally, machine learning approaches treat the changepoint detection problem through either a regression or classification \cite{Kotsiantis2007} function approximation and seek to identify the category an observation $\boldsymbol{z}_t$ belongs to. In particular, neural network models have recently been used for the purpose of sequence classification\cite{Ahad2021} and time series anomaly detection\cite{Pang2020, Munir2018}. In the online setting, one may flag a changepoint when a trained classifier or cluster predicts a different category of $\boldsymbol{z}_t$ from $\boldsymbol{z}_{t-1}$. 

\section{Dataset}\label{sec:dataset}

\subsection{Kolmogorov Flow}

We utilize the two-dimensional Navier-Stokes equations with Kolmogorov forcing (commonly denoted Kolmogorov flow) as given by
\begin{align}
    \frac{\partial \mathbf{u}}{\partial t} + \mathbf{u}. \nabla \mathbf{u} + \Delta p &= \frac{1}{\text{Re}} \nabla^2 \mathbf{u} + \sin{(ny)} \hat{\mathbf{x}} \nonumber \\ 
    \nabla . \mathbf{u} &= 0
\end{align}
where $\mathbf{u} = [u,v]$ is the velocity vector, $p$ is the pressure, $n$ is the wavenumber of the forcing and $\hat{\mathbf{x}}$ is the unit vector in the $x$ direction. Here $Re = \frac{\sqrt{\chi}}{\nu} \frac{L_y}{2 \pi}^{3/2}$, where $\chi$ is the dimensional forcing amplitude, $\nu$ is the kinematic viscosity, and $L_y$ is the size of the domain in the $y$ direction. We consider a square periodic domain for this problem with a side length of $2 \pi$. Our raw state for machine learning computation and analyses in this paper is the vorticity given by $\omega = \nabla \times \mathbf{u}$. An important quantity of interest for this system is the total kinetic energy given by
\begin{align}
    KE = \frac{1}{2} \langle \mathbf{u}^2 \rangle_V
\end{align}
where the subscript $V$ corresponds to the volume average. We use the data and numerical methods presented in \cite{de2023data} for the purpose of our experiments and analyses. In particular we use a Reynolds number of 14.4 to obtain trajectories of data where there are intermittent bursting events. We obtain trajectories of the state and for bursting events from De Jes\'{u}s and Graham\cite{de2023data}.


\section{Methods}\label{sec:methods}

This section is devoted to the exposition of our proposed numerical techniques for scalable changepoint detection of high-dimensional systems. Our two methods follow different intuitions to accomplish the same goal. In the first, we assume that the observed dynamics are generated by some underlying generative process and detect sudden changes in the parameterizations of this process in an online manner. In the second, we assume that a drastic change in incoming information can be characterized as a sudden drop in the predictability of a predictive model that is observed a sliding window of inputs. We proceed with the mathematical descriptions of the different techniques below.

\subsection{Bayesian online changepoint detection}

Our first approach to changepoint detection is given by Bayesian online changepoint detection\cite{BOCPD} (BOCD). BOCD models changepoints in time-series data as sudden shifts in the underlying distribution that may potential indicate transitions between different physical regimes in the context of dynamical systems. Moreover, BOCD is able to perform this detection task in an online manner as shall be discussed in the following. BOCD assumes that a sequence of observations can be partitioned into non-overlapping states, where each state $\rho$ is characterized by an independent and identically distributed (i.i.d.) set of observations $\boldsymbol{z}_t$, drawn from a probability distribution $P(\boldsymbol{z}_t|\eta_\rho)$ with i.i.d. parameters $\eta_\rho,\ \rho=1,2,...$. To facilitate change point detection, BOCD introduces an auxiliary variable $r_t$, which represents the elapsed time since the last change point. Given the run length $r_t$ at time $t$, the run length at the next time point can either be $0$, indicating the occurrence of a change point, or increase by $1$, if the current state persists. The target then becomes to find the posterior probability:
\begin{eqnarray}
P(r_t|\boldsymbol{z}_{1:t}) =\frac{P(r_t,\boldsymbol{z}_{1:t})}{P(\boldsymbol{z}_{1:t})}
\end{eqnarray}
where
\begin{eqnarray}
P(r_t,\boldsymbol{z}_{1:t}) = \sum_{r_{t-1}}P(r_t |r_{t-1} )P(z_{t} |r_{t-1} ,\boldsymbol{z}_t^{(r)})P(r_{t-1} ,\boldsymbol{z}_{1:t-1} ).
\end{eqnarray}
Here, $\boldsymbol{z}_t^{(r)}$ indicates the set of observations associated with the run $r_t$. $P(r_t |r_{t-1} )$ is the prior, $P(z_{t} |r_{t-1} ,\boldsymbol{z}_{1:t})$ is the predictive distribution which only depend on the recent data $\boldsymbol{z}_t^{(r)}$, and $P(r_{t-1} ,\boldsymbol{z}_{1:t-1})$ are and recursive components of the equation. The conditional prior is nonzero at only two outcomes ($r_t = 0$ or $r_t = r_{t-1} + 1$).
\begin{equation}
  P(r_t |r_{t-1} ) =
    \begin{cases}
      H(r_{t-1}+1) & \text{if $r_t=0$}\\
      1-H(r_{t-1}+1) & \text{if $r_t=r_{t-1}+1$}\\
      0 & \text{otherwise}.
    \end{cases}       
\end{equation}
The function $H(\tau)$ is hazard function. $H(\tau) = \frac{P(g=\tau)}{\sum_{t=\tau}^\infty P(g=t)}$. We assume that the process is Markovian and configure $P(g)$ as a discrete exponential distribution with given parameter $\lambda$. In addition, the hazard function is held constant at $H(\tau)=\frac{1}{\lambda}$.

The predictive distribution $P(z_{t} |r_{t-1} ,\boldsymbol{z}_t^{(r)})$ may become exceptionally intractable to calculate in the generic case, and therefore, a conjugate exponential model is used for most computationally efficiency. Exponential family likelihoods allow inference with a finite number of sufficient statistics $\eta$ which can be updated incrementally as data arrives. The fast update of the predictive distribution by using a conjugate model allows for fast online computation as new data is observed. A recursive message passing algorithm can then be implemented to calculate $P(r_t|\boldsymbol{z}_{1:t})$. Our computation is initialized with $P(r_0=0)=1$, and carefully selected hyperparameters for the predictive distribution. Several schemes can be adopted to extract the location of change points\cite{ocp}:
\begin{itemize}
    \item Threshold based scheme: When new data arrives at time $t$, the posterior distribution $P(r_t|\boldsymbol{z}_{1:t})$ is evaluated. If the posterior probability $P(r_t=k|\boldsymbol{z}_{1:t})$ exceeds a given threshold, we flag a changepoint at time $t-k$.
    \item Maximum probability based scheme: Similar to the threshold-based scheme. When new data arrives, we calculate the posterior distribution of $r_t$ and find the maximum probability $P(r_t|\boldsymbol{z}_{1:t})$ at $r_t= k$, we flag a changepoint at time $t-k$.
    \item Most probable changepoints set based scheme: While calculating the posterior distribution $P(r_t|\boldsymbol{z}_{i:t})$, we track and record the probability of the most likely set of changepoints $MP(t)$. In other words, if the maximum of $AP(r_t)=P(r_t|\boldsymbol{z}_{i:t})+MP(t)$ is found at $r_t=k$, we flag a changepoint at time $t-k$.
\end{itemize}

In the above description, we have assumed a multivariate observation for $\boldsymbol{z}$ which implies that a covariance matrix update must be performed at each timestep, within the framework of the conjugate model update. We will demonstrate how an intelligent dimension reduction technique of our original state data can further avoid this computational expense. We will demonstrate how a `disentangled' representation of the state vector, where individual dimensions of the state are independent of each other, can be used to compute the multivariate prediction distribution while avoiding the need to update the covariance matrix. The prediction log-likelihood, as well as parameter updates, can then be implemented in parallel for every dimension. The sum of prediction log-likelihoods can then be used to calculate the posterior distribution of the run length, and changepoint detection can be carried out using one of the three schemes mentioned above. In our experiments, we use the third approach for the purpose of robustness.

\begin{algorithm}[H]
\caption{Bayesian online changepoint detection}
\begin{algorithmic}[1]
 \State $P(r_0=0)=1,\ MP(0) = [0,0]$ 
 \State $\mu_0 = \mu_0,\ \kappa_0=\kappa_0,\ \alpha_0=\alpha_0,\ \beta_0=\beta_0,\ \lambda=\lambda$
 \For {Each $t$}
	\For {$i=1,\ldots,d$}
		\State Calculate predictive probability: $\pi(\boldsymbol{z}_t(i))=P(\boldsymbol{z}_t(i)|\boldsymbol{z}_{t-1}^{(r)}(i),r_{t-1})$
		\State Update sufficient statistics according to Murphy \cite{Murphy2007}
	\EndFor	
        \State Calculate predictive probability: $\pi(\boldsymbol{z}_t)=\prod^d_i P(\boldsymbol{z}_t(i)|\boldsymbol{z}_{t-1}^{(r)}(i),r_{t-1})$
	\State Calculate growth probability: $P(r_t=r_{t-1}+1,\boldsymbol{z}_{1:t})=P(r_{t-1},\boldsymbol{z}_{1:t-1})\pi(\boldsymbol{z}_t)(1-H(r_{t-1}))$
        \State Calculate changepoint probability: $P(r_t=0,\boldsymbol{z}_{1:t})=\sum_{r_{t-1}}P(r_{t-1},\boldsymbol{z}_{1:t-1})\pi(\boldsymbol{z}_t)H(r_{t-1})$
	\State Normalize probability of $r_t$: $P(r_t|\boldsymbol{z}_{1:t}) = \frac{P(r_t,\boldsymbol{z}_{1:t})}{\sum_{r_t}P(r_t,\boldsymbol{z}_{1:t})}$
        \State Adjust probability of $r_t$: $AP(r_t)=P(r_t|\boldsymbol{z}_{1:t}+MP(t)$
        \State Update probability of most probable set of changepoints: $MP = [max(AP(r_t)),MP]$
        \If{$argmax_{r_t} AP(r_t)=0$}
          \State flag a changepoint at time $t$
        \EndIf
 \EndFor
\end{algorithmic} 
\end{algorithm}

\subsection{Supervised learning based changepoint detection}

In this section, we highlight our other approach to changepoint detection that relies on supervised learning techniques for time-series forecasting. Here, the intuition is to use the (lack of) model predictability to flag potential changes in the underlying regime of the observed dynamics. In this work, we rely on the use of the well-known Long Short-term Memory neural network (LSTM)\cite{Hochreiter1997},\cite{Gers2002},\cite{Gers2000} for our supervised learning model. LSTMs have been sccessfully applied to time series forecasting. It is a type of Recurrent Neural Network (RNN)\cite{Werbos1990},\cite{Williams1989} that can model sequential data with long-term dependencies. Unlike traditional RNNs, which suffer from the vanishing or exploding gradient problem when trained on long sequences\cite{Mozer1992}, LSTMs use a gating mechanism to selectively remember or forget information at each time step.

The basic building block of an LSTM cell consists of three gates (input, forget, and output) and a memory cell. The input gate decides which information to add to the memory cell, the forget gate decides which information to discard from the memory cell, and the output gate decides which information to output from the memory cell. These gates are controlled by sigmoid activation functions that output values between 0 and 1, which act as weights for the input, forget, and output operations.

Mathematically, an LSTM cell can be defined as follows:
\begin{align*}
\boldsymbol{i}_t &= \sigma(\boldsymbol{W}_{ii}\boldsymbol{z}_t + \boldsymbol{b}_{ii}+ \boldsymbol{W}_{hi}\boldsymbol{h}_{t-1} + \boldsymbol{b}_{hi}) \\
\boldsymbol{f}_t &= \sigma(\boldsymbol{W}_{if}\boldsymbol{z}_t + \boldsymbol{b}_{if}+ \boldsymbol{W}_{hf}\boldsymbol{h}_{t-1} + \boldsymbol{b}_{hf}) \\
\boldsymbol{o}_t &= \sigma(\boldsymbol{W}_{io}\boldsymbol{z}_t +\boldsymbol{b}_{io}+ \boldsymbol{W}_{ho}\boldsymbol{h}_{t-1} + \boldsymbol{b}_{ho}) \\
\boldsymbol{g}_t &= \tanh(\boldsymbol{W}_{ig}\boldsymbol{z}_t + \boldsymbol{b}_{ig}+ \boldsymbol{W}_{hg}\boldsymbol{h}_{t-1} + \boldsymbol{b}_{hg}) \\
\boldsymbol{c}_t &= \boldsymbol{f}_t\odot\boldsymbol{c}_{t-1} + \boldsymbol{i}_t\odot\boldsymbol{g}_t \\
\boldsymbol{h}_t &= \boldsymbol{o}_t\odot\tanh(\boldsymbol{c}_t)
\end{align*}

where $\boldsymbol{z}_t$ is the input at time step $t$, $\boldsymbol{h}_t$ is the hidden state at time step $t$, $\boldsymbol{c}_t$ is the memory cell at time step $t$, $\boldsymbol{i}_t$, $\boldsymbol{f}_t$, and $\boldsymbol{o}_t$ are the input, forget, and output gates at time step $t$, respectively, and $\boldsymbol{g}_t$ is the candidate memory at time step $t$. $\sigma$ is the sigmoid activation function, $\odot$ denotes element-wise multiplication, and $\boldsymbol{W}$ and $\boldsymbol{b}$ are the learnable weight matrices and bias vectors, respectively.

Given a time-series with $k$ steps $\boldsymbol{z}_{history} = (\boldsymbol{z}_t, \boldsymbol{z}_{t-1},...,\boldsymbol{z_{t-k+1}})$, the objective is to predict the next $k$ steps $\boldsymbol{z}_{pred} = (\boldsymbol{z}_{t+1},...,\boldsymbol{z_{t+k}})$ using an LSTM model. The LSTM model is trained on a set of sequences of fixed length $k$ by minimizing the L1-norm, which provides model robustness against anomalies, between the predicted and actual values. This allows for a model that predicts closer to the median of the trajectory and predicting sudden changes in state trajectories are not emphasized during the prediction process. Once the LSTM model is trained, a sequence of predicted values $\hat{\boldsymbol{z}}_{t+1},...,\hat{\boldsymbol{z}}_{t+k}$ is generated. The maximum value of L2 norm of prediction error in every step is used to compute the anomaly score as follows:
\begin{align}
    \label{anom_score}
    AS = \max_{j=1,...,k}{||\hat{\boldsymbol{z}_j} - \boldsymbol{z}_j||_2}
\end{align}

To detect change points, the anomaly score is calculated when the $k$th step arrives. If the anomaly score exceeds a threshold value, the transition is assumed to begin within the next $k$ steps. During the transition phase, the prediction error will be constantly high.

\subsection{Dimension Reduction with $\beta$-TCVAE}

In this section, we introduce a dimensionality reduction technique to accelerate the changepoint detection techniques introduced previously, particularly for BOCD. We begin by noting that performing changepoint detection by observing the high-dimensional state obtained by numerically solving the Navier-Stokes equations introduced for Kolmogorov flow are infeasible from the perspective of scaling to practical tasks. For BOCD, one needs to model and update a covariance matrix for every new snapshot of the state and a high-dimensional observation with correlated state components significantly complicates this task computationally. Therefore, we seek to address the correlation between components while reducing the dimensionality of the problem by applying dimension reduction and disentanglement techniques in machine learning. 

Specifically, we utilize the $\beta$-TCVAE\cite{Chen2019} approach which is a state-of-the-art technique to jointly accomplish this task. The $\beta$-TCVAE is a derivative of a variational autoencoder (\textbf{VAE})\cite{Higgins2017, Kingma2013}, which is a latent (reduced dimensionality) variable model. It consists of a inference network and a generative network. The inference network maps observed high-dimensional data variables $\omega_t\in \mathbb{R}^m$ to latent low dimensional variables $\boldsymbol{z}_t\in\mathbb{R}^k$, where $k<<m$. This inference network is formulated as $q(\boldsymbol{z}_t|\omega_t)$. The generative network generates variables from a latent distribution, namely $p(\omega_t|\boldsymbol{z}_t)$. A VAE-type model pairs the generator with inference network and training is done by optimizing the tractable averaged evidence lower bound (ELBO) over the empirical distribution (with $\beta = 1$):
\begin{align}
    \mathcal{L}_{\beta} = \frac{1}{T}\sum^T_{t=1}(\mathbb{E}_q[\text{log}p(\omega_t|\boldsymbol{z})]-\beta \textbf{KL}(q(\boldsymbol{z}|\omega_t)||p(\boldsymbol{z}))).
    \label{elbo}
\end{align}
In the \textbf{$\beta$-TCVAE}, we define $q(\boldsymbol{z},\omega_t) = q(\boldsymbol{z}|\omega_t)p(\omega_t) = q(\boldsymbol{z}|\omega_t)\frac{1}{T}$ and take $q(\boldsymbol{z}) = \sum_{t=1}^T q(\boldsymbol{z}|\omega_t)p(\omega_t)$ as the aggregated posterior. This decomposes the KL term in \eqref{elbo} into:
\begin{eqnarray}
\mathbb{E}_{p(\omega_t)}[\textbf{KL}(q(\boldsymbol{z}|\omega_t)||p(\boldsymbol{z}))] &&= \textbf{KL}(q(\boldsymbol{z},\omega_t)||q(\boldsymbol{z})p(\omega_t)) \nonumber\\
&&+ \textbf{KL}(q(\boldsymbol{z})||\prod_j q(\boldsymbol{z}_j)) \nonumber\\
&&+ \sum_j\textbf{KL}(q(\boldsymbol{z}_j)||p(\boldsymbol{z}_j)),
\end{eqnarray}
and provides a new ELBO:
\begin{align}
    \mathcal{L}_{\beta-TC} &&= \mathbb{E}_q[\text{log}p(\omega_t|\boldsymbol{z})]-\alpha I_q(\boldsymbol{z},\omega_t)- \beta\textbf{KL}(q(\boldsymbol{z})||\prod_j q(\boldsymbol{z}_j)) \nonumber\\
    &&- \gamma\sum_j\textbf{KL}(q(\boldsymbol{z}_j)||p(\boldsymbol{z}_j))
    \label{tcvae}
\end{align}
where $I_q(\boldsymbol{z},\omega_t)=\textbf{KL}(q(\boldsymbol{z},\omega_t)||q(\boldsymbol{z})p(\omega_t))$ is the mutual information between $\boldsymbol{z}$ and $\omega_t$. Chen\cite{Chen2019} has verified that better trade-off between ELBO and disentanglement can be realized by tuning hyperparameter $\beta$ (set $\alpha=\gamma=1$).

In order to ensure that losses due to inaccurate reconstruction were weighted fairly across the different variables, the data is scaled to have a zero mean and unit variance. To obtain a low-dimensional representation of the solution field, we sue a two-dimensional convolutional framework with multiple strided filters. To reduce the dimensionality of the input image to a size of $4\times 2$ degrees of freedom in the latent space, fully connected convolutional layers are utilized. These degrees of freedom represent the mean $\mu$ and standard deviation $\sigma$ of the latent space, which is assumed to be a multivariate normal distribution. A 4-dimensional sample is then sampled from the latent distribution using the reparameterization trick. Following this, the 4-dimensional state is convolved and upsampled several times to return to the dimensionality of the full-order field. Each layer consisted of rectified linear (ReLU) activations and utilizes zero-padding at the edges of the domain for the purpose of convolution. Our network is trained using the loss function in equation \eqref{tcvae} with a batch size of $1024$. We used the Adam optimizer with a learning rate of 0.001. In this case, we assumed that the latent variable $\boldsymbol{z}$ and the posterior $p(\omega_t|\boldsymbol{z}_t)$ followed a normal distribution. The cross entropy $\mathbb{E}_q[\text{log}p(\omega_t|z)]$ was simplified as the reconstruction error, which was represented using standard Mean Squared Error (MSE) between the input and output.  Each convolutional layer in the autoencoder utilizes a ReLU activation function, with the exception of the output layer and the final layer of the encoder. We use 49991 snapshots (first half snapshots) for training and tested the trained $\beta-TCVAE$ model on the remaining snapshots. Specific details of the architecture, such as the number of channels in each layer and functions used to calculate the loss function, can be found in the supporting code.

\section{Hyperparameter tuning}\label{sec:hyperparam}

As is common in most scientific machine learning studies, extensive hyperparameter tuning is necessary for the appropriate deployment of the previously discussed models. In this section, we outline these tuning investigations. 

\subsection{Hyperparameter tuning for BOCD}

In our study, we assume data $\boldsymbol{z}$ are normally distributed with unknown variance $\sigma^2$ and mean $\mu$. We use a conjugate prior distribution, the Normal-Gamma distribution with parameters $\mu_0,\ \kappa_0,\ \alpha_0,$ and $\beta_0$. The unknown $\mu$ follows a normal distribution with mean $\mu_0$ and variance $\sigma^2/\kappa_0$, while the unknown $\sigma^{-2}$ follows gamma distribution with shape parameter $\alpha_0$ and rate parameter $\beta_0$. Subsequently, the posterior predictive distribution $P(\boldsymbol{z}_t|\boldsymbol{z}_t^{(r)})$ is a Student's-T distribution Murphy (2007)\cite{Murphy2007}, and the parameters are updated accordingly with new data. As the data is standardized, we set the prior mean set to $\mu_0 = 0$ and $\kappa_0=1$. Also, we remind the reader that the hazard function is constant at $H(\tau)=\frac{1}{\lambda}$.

To investigate the impact of hyperparameters on the performance of BOCD, we varied three hyperparameters: $\alpha_0$, $\beta_0$, and $\lambda$. Here, $\lambda$ indicates the expected interval between two successive changepoints and was selected from the set of values $\{5,10,50,100\}$. The other parameters, $\alpha_0$ and $\beta_0$ are varied across the sets of values $\{0.01, 0.1, 1, 10\}$ which represent possible values of unknown $\sigma^{-2}$. To determine the optimal configuration, we performed a grid search, where we systematically tested all possible hyperparameter combinations. We note that more intelligent hyperparameter searches can be performed to optimize the discovery of parameters\cite{Turner2009}.


To evaluate BOCD performance with different hyperparameter settings, we use the $F-$score metric\cite{Rijsbergen1979}. This metric is commonly used in change point detection studies and is a special form of evaluation introduced in [\cite{killick12},\cite{vandenBurg2022},\cite{Truong2020}]. The set of change points is denoted by the ordered set $\mathcal{T} = \{\tau_1,...,\tau_{n}\}$ with $\tau_i \in [1,T]$ for $i = 1,...,n$ and $\tau_i < \tau_j$ for $i < j$. Then $\mathcal{T}$ implies a partition, $\mathcal{G}$, of the interval $[1, T ]$ into disjoint sets $\mathcal{A}_j$, where $\mathcal{A}_j$ is the segment from $\tau_{j-1}$ to $\tau_j-1$ for $j = 1,...,n+1$. For notational convenience, we set $\tau_0 = 1$ and $\tau_{n+1} = T +1$.

The $F-$score measures the performance of change point detection as a classification problem between `change point' and `non-change point' classes. To allow for minor discrepancies, we define a margin of error $M\geq0$ around the true change point location. For the set of ground truth locations $\mathcal{T}$ provided by experts and change point locations $\mathcal{T'}$ detected by the algorithm, we define the set of true positives $TP(\mathcal{T} , \mathcal{T}')$ to be those $\tau \in \mathcal{T}$ for which $\exists \tau'\in\mathcal{T}'$ such that $|\tau-\tau'|\leq M$, while ensuring that only one $\tau'\in\mathcal{T'}$ can be used for a single $\tau\in\mathcal{T}$. Furthermore, we define precision and recall as:
\begin{eqnarray}
&&P = \frac{|TP(\mathcal{T}, \mathcal{T}')|}{|\mathcal{T}'|},\\
    &&R=\frac{|TP(\mathcal{T}, \mathcal{T}')|}{|\mathcal{T}|},\\
    &&F = \frac{2PR}{P+R}
\end{eqnarray}
In the experiments we use a margin of error of $M = 5$ points.

\subsection{Hyperparameter tuning for supervised learning based changepoint detection}

For the model predictability based changepoint detection, we identify two hyperparameters for tuning. These are the prediction window size $k$, and $thr$, the threshold for flagging a changepoint using an anomaly score. We choose $k = 5, 10, 15$ in our experiments. While the threshold is determined after we an overall examination of anomaly scores predicted by the trained models.

\subsection{Hyperparameter tuning for $\beta$-TCVAE}

The $\beta$-TCVAE based dimensionality reduction technique is used to obtain a trade-off between reconstruction error and disentanglement. We remind the reader that the latter is important for the accelerated detection of changepoints in BOCD. Moreover, an information latent space, given by low reconstruction error, is important for accurate capture of variations of the high-dimensional state that may represent changepoints. Therefore, the parameter $\beta$, which controls this trade-off requires careful selection. To this end, we trained models with $\beta$ values of 1, 5, 10, and 20 on the training set, and calculated the model's reconstruction error on the validation set. 

\begin{figure*}[!htb]
 \includegraphics[width=1\textwidth]{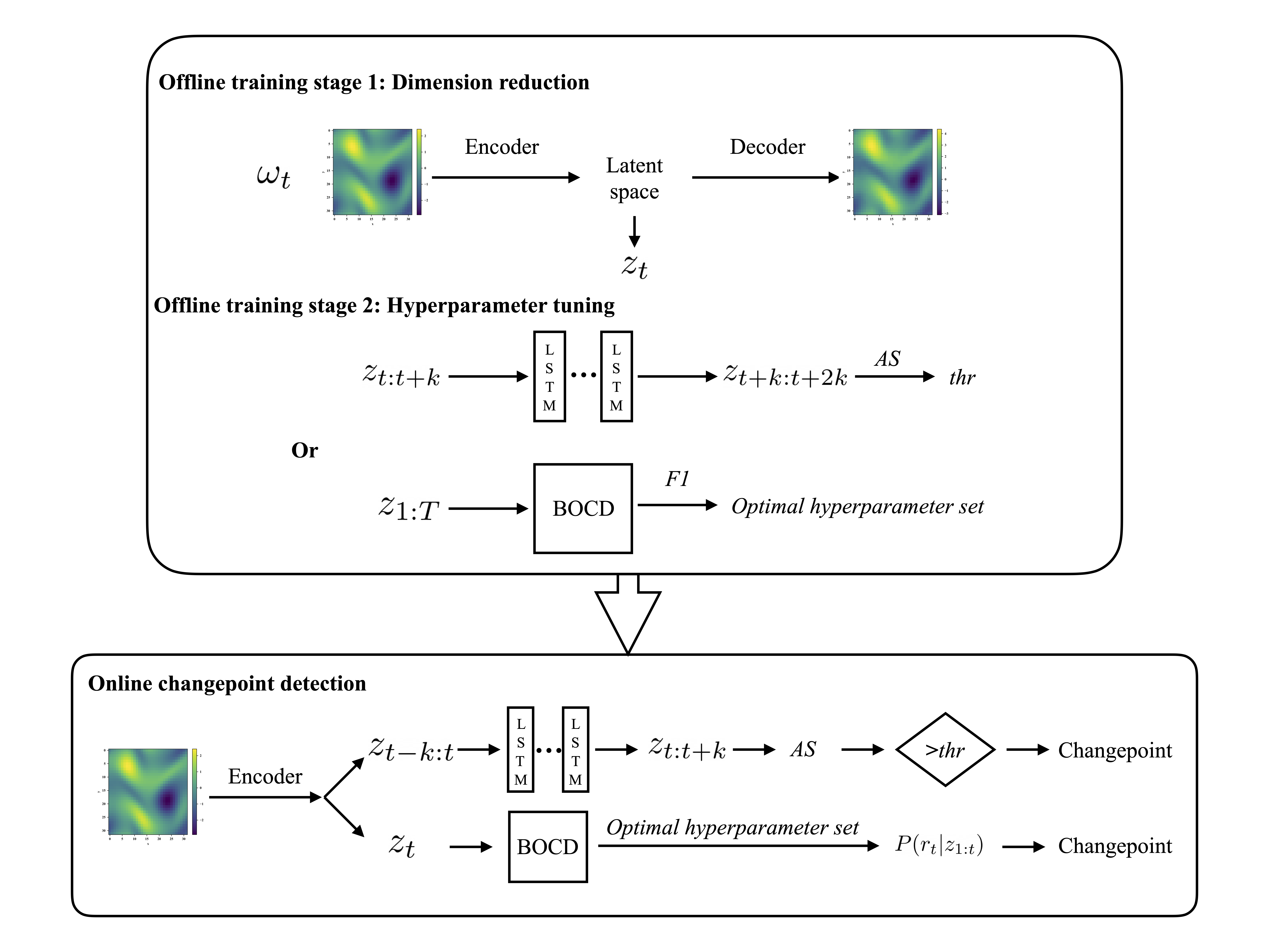}
 \caption{\label{workflow}Workflow for the entire algorithm: Vorticity snapshots are compressed using a variational autoencoder to obtain low-dimensional representations. The latter are used, after optimal selection of hyperparameters, to flag anomalies in real-time.}
 \label{fig:workflow}
\end{figure*}

\section{Results}\label{sec:results}
 
We outline results from our studies in the following.

\subsection{Dimensionality reduction and disentanglement}

\begin{figure*}[!htb]
    \centering
        \mbox{
        \subfigure[Latent feature $z_1$ and kinetic energy.]{\includegraphics[width=0.4\textwidth]{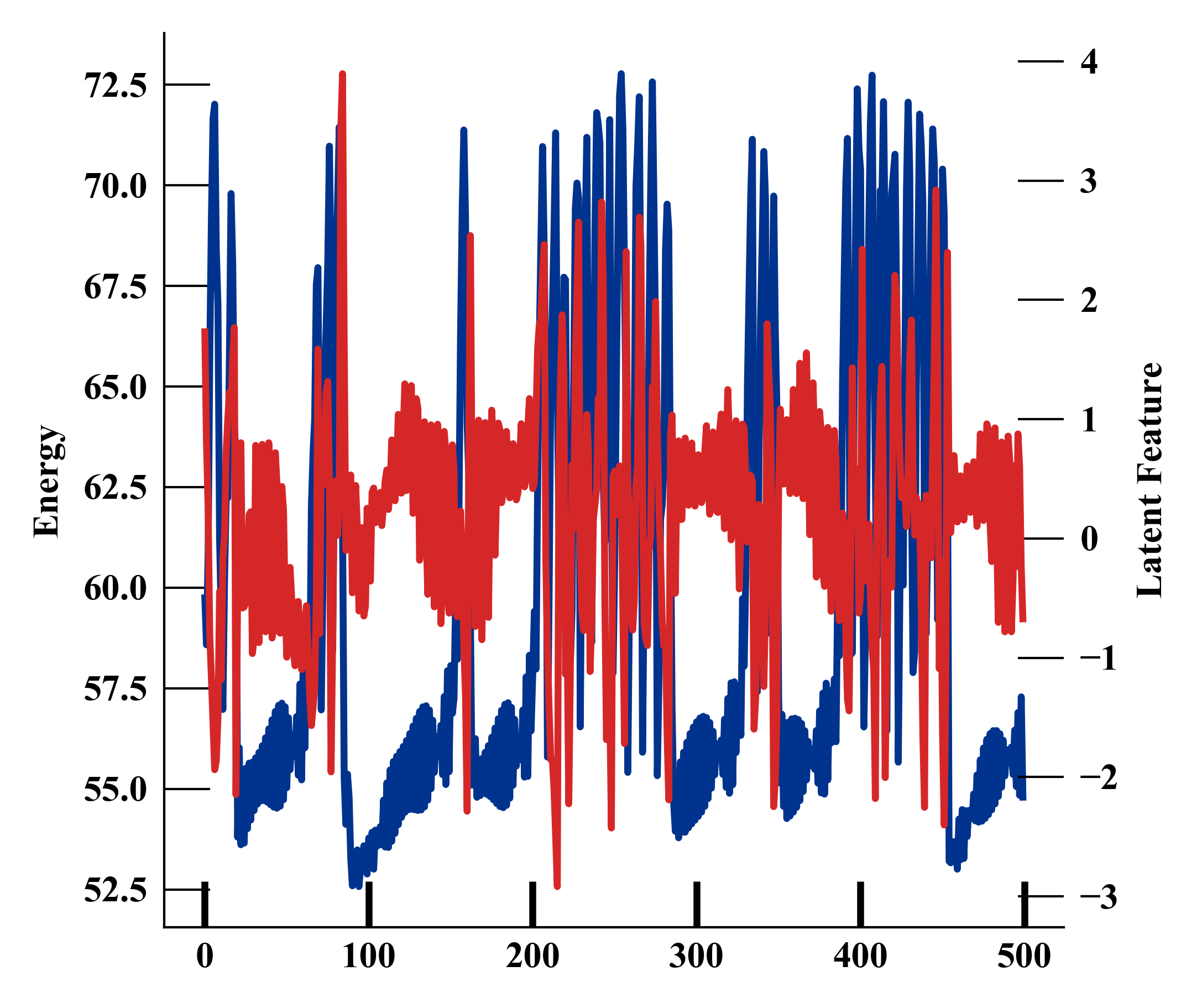}}
        \subfigure[Latent feature $z_2$ and kinetic energy.]{\includegraphics[width=0.4\textwidth]{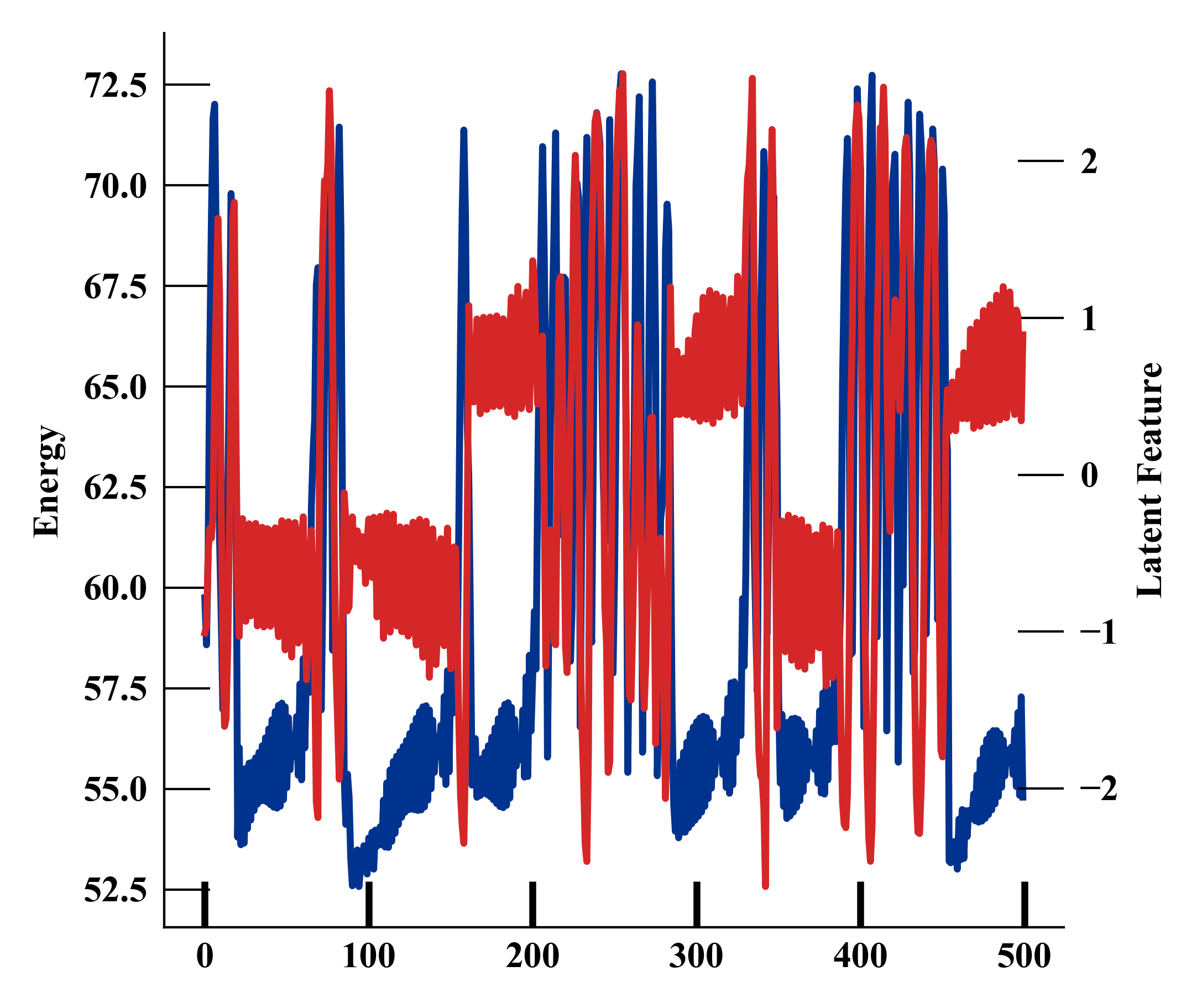}}
        } \\
        \mbox{
        \subfigure[Latent feature $z_3$ and kinetic energy.]{\includegraphics[width=0.4\textwidth]{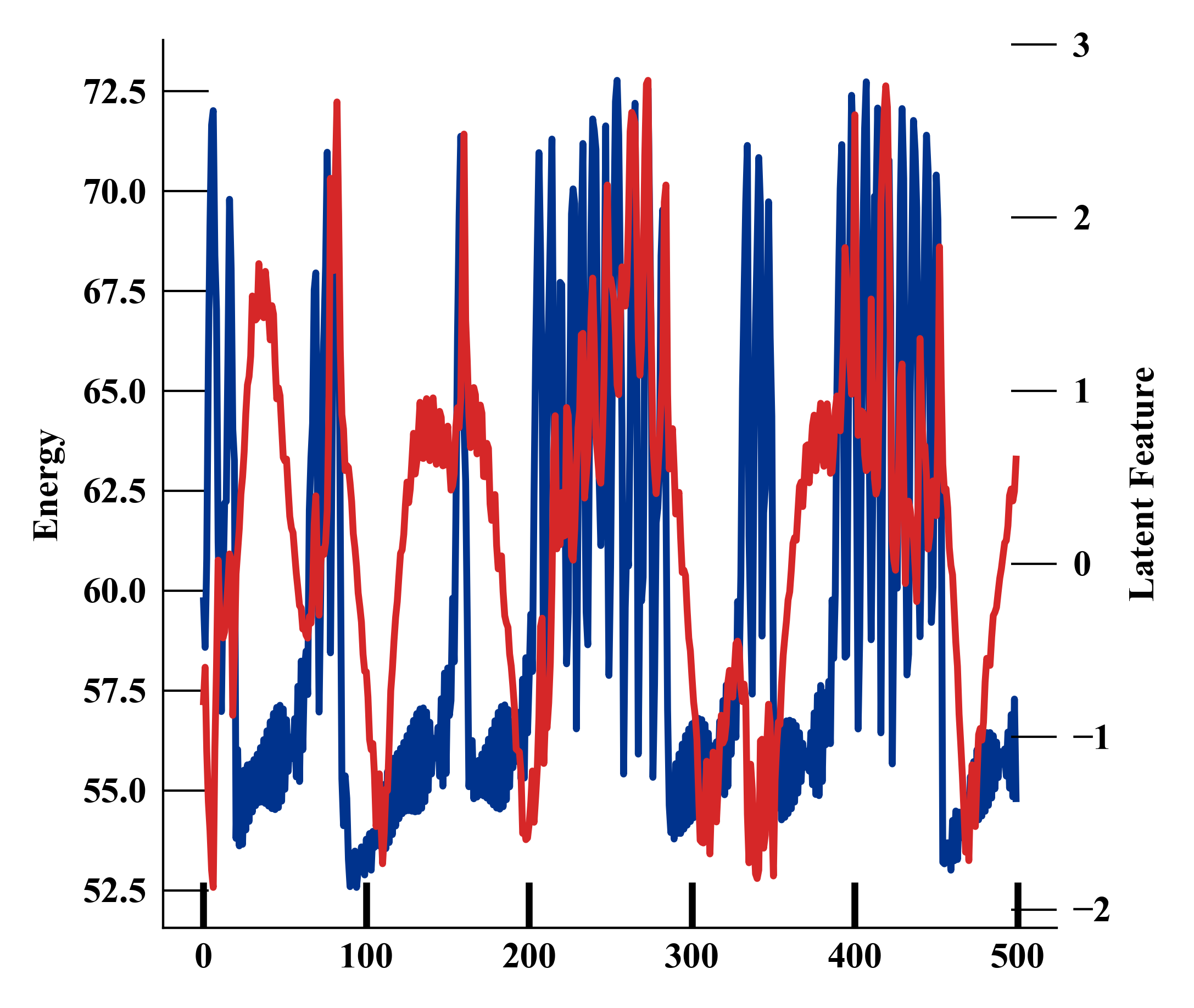}}
        \subfigure[Latent feature $z_4$ and kinetic energy.]{\includegraphics[width=0.4\textwidth]{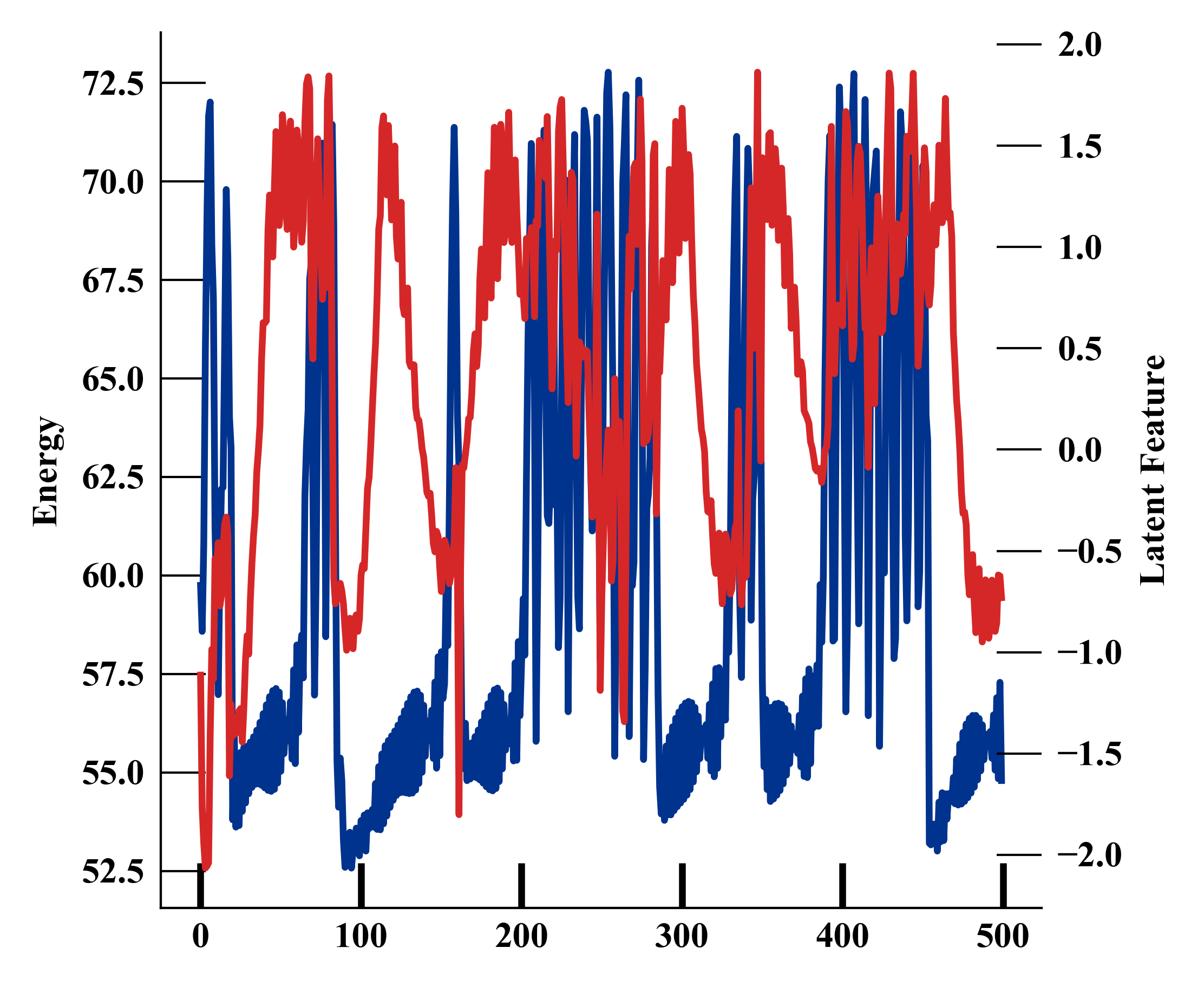}}
        }
    \caption{Relation between the latent features and kinetic energy. Curves for kinetic energy in each subplot are in blue.}
    \label{fig:energy_latent}
\end{figure*}

First, we outline our studies with successfully compressing the high-dimensional state, given by the vorticity of the Kolmogorov flow $\omega$, using the $\beta$-TCVAE. The snapshots were successfully compressed into low-dimensional representations, and the relationship between the latent dimension and kinetic energy is shown in Fig.\ref{fig:energy_latent}. It may clearly be observed that the latent dimensions, though disentangled, are highly correlated with transitions between regular and intermittent bursting events.

\begin{figure*}[!htb]
    \centering
    \mbox{
    \subfigure[Original]{
    \includegraphics[width=0.4\textwidth]{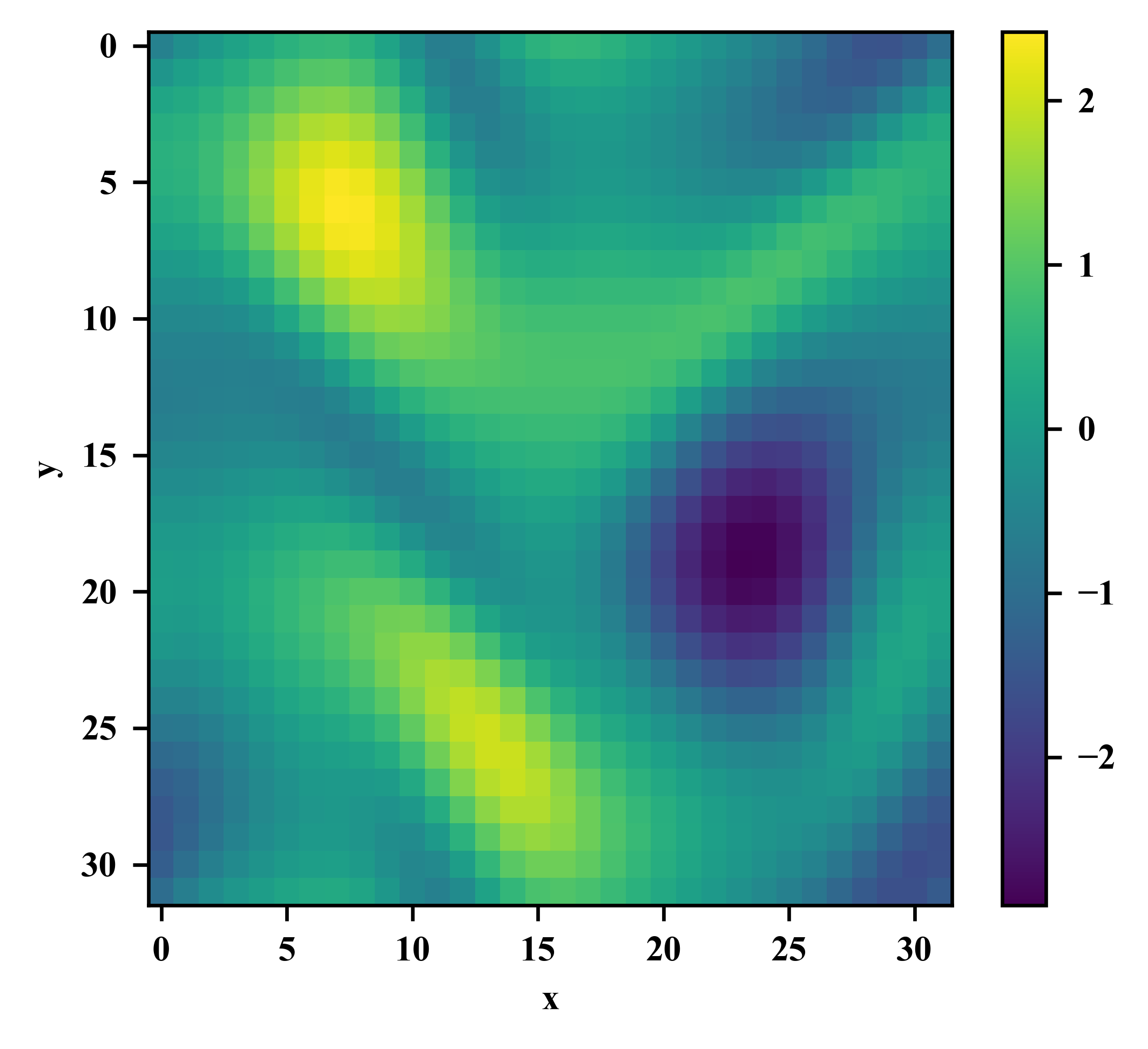}
    \label{fig:original}
    }
    \subfigure[Reconstructed]{\includegraphics[width=0.4\textwidth]{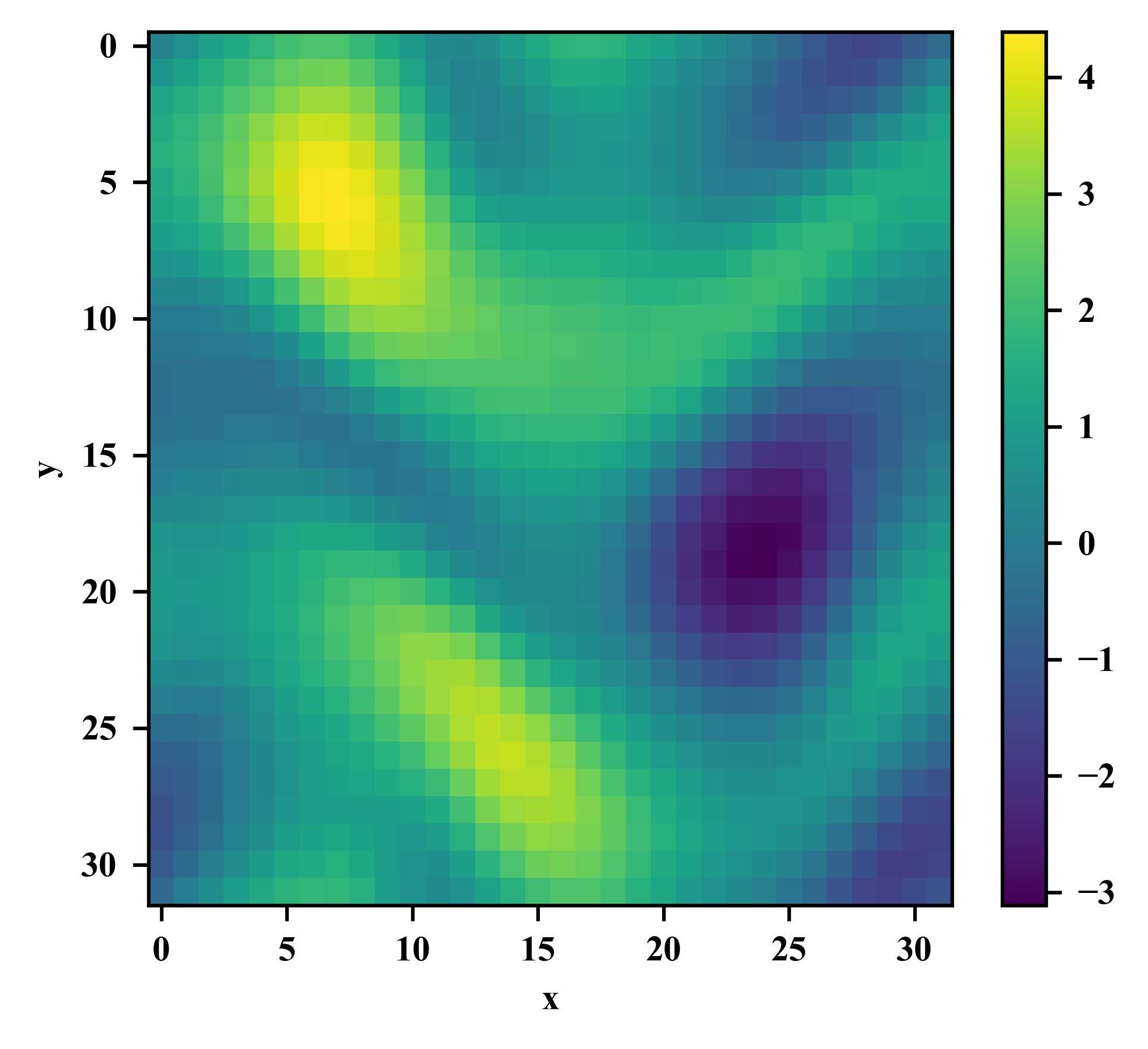}
    \label{fig:reconstructed}
    }
    }
    \caption{Comparison between origin snapshot and reconstructed snapshots of vorticity obtained from $\beta-$TCVAE training.}
    \label{fig:reconstruction_vae}
\end{figure*}

\begin{figure*}[!htb]
    \centering
    \mbox{
    \subfigure[Reconstructed snapshots with perturbations in $z_1$]{
    \includegraphics[width=0.4\textwidth]{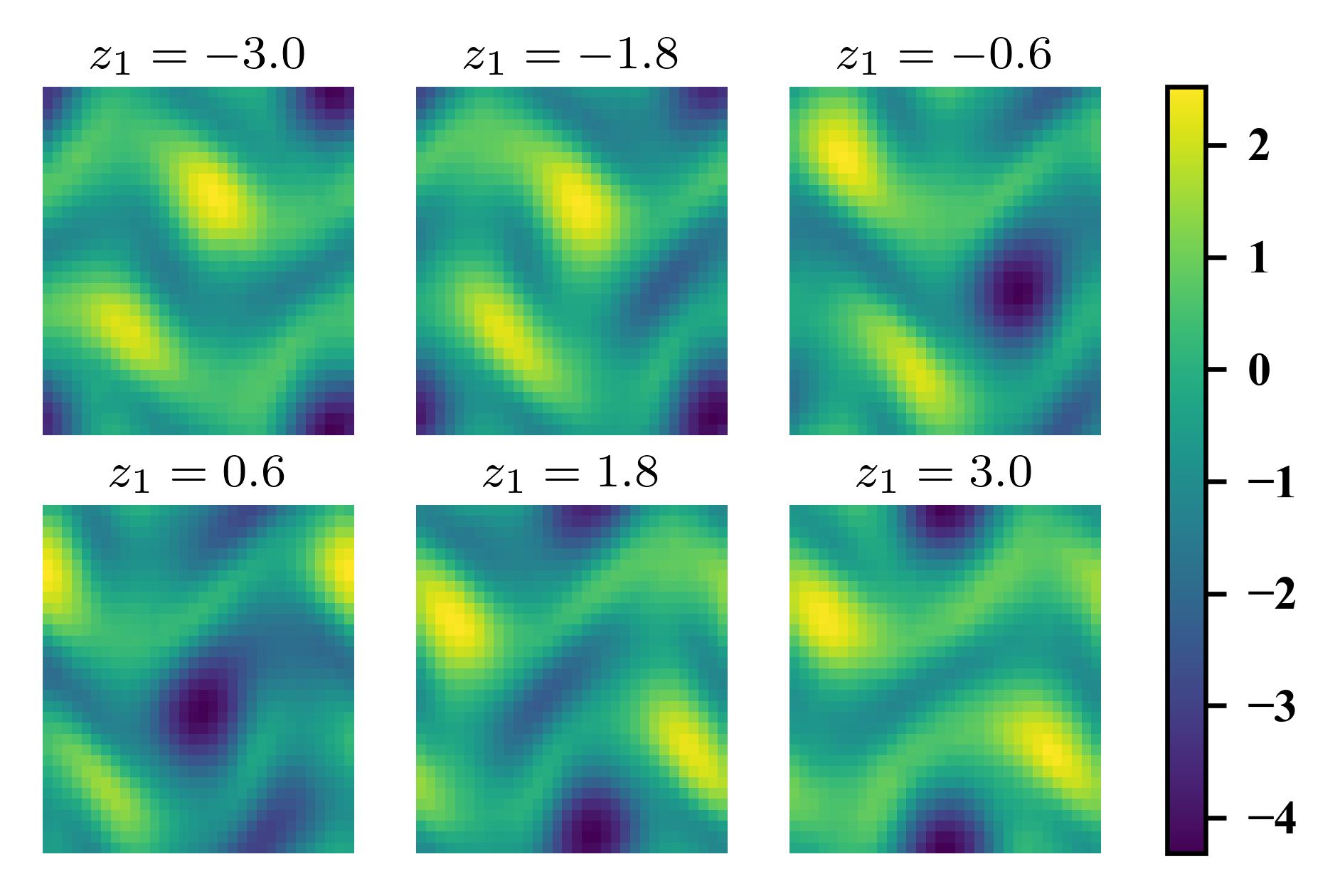}
    \label{fig:z1}
    }
    \subfigure[Reconstructed snapshots with perturbations in $z_2$]{
    \includegraphics[width=0.4\textwidth]{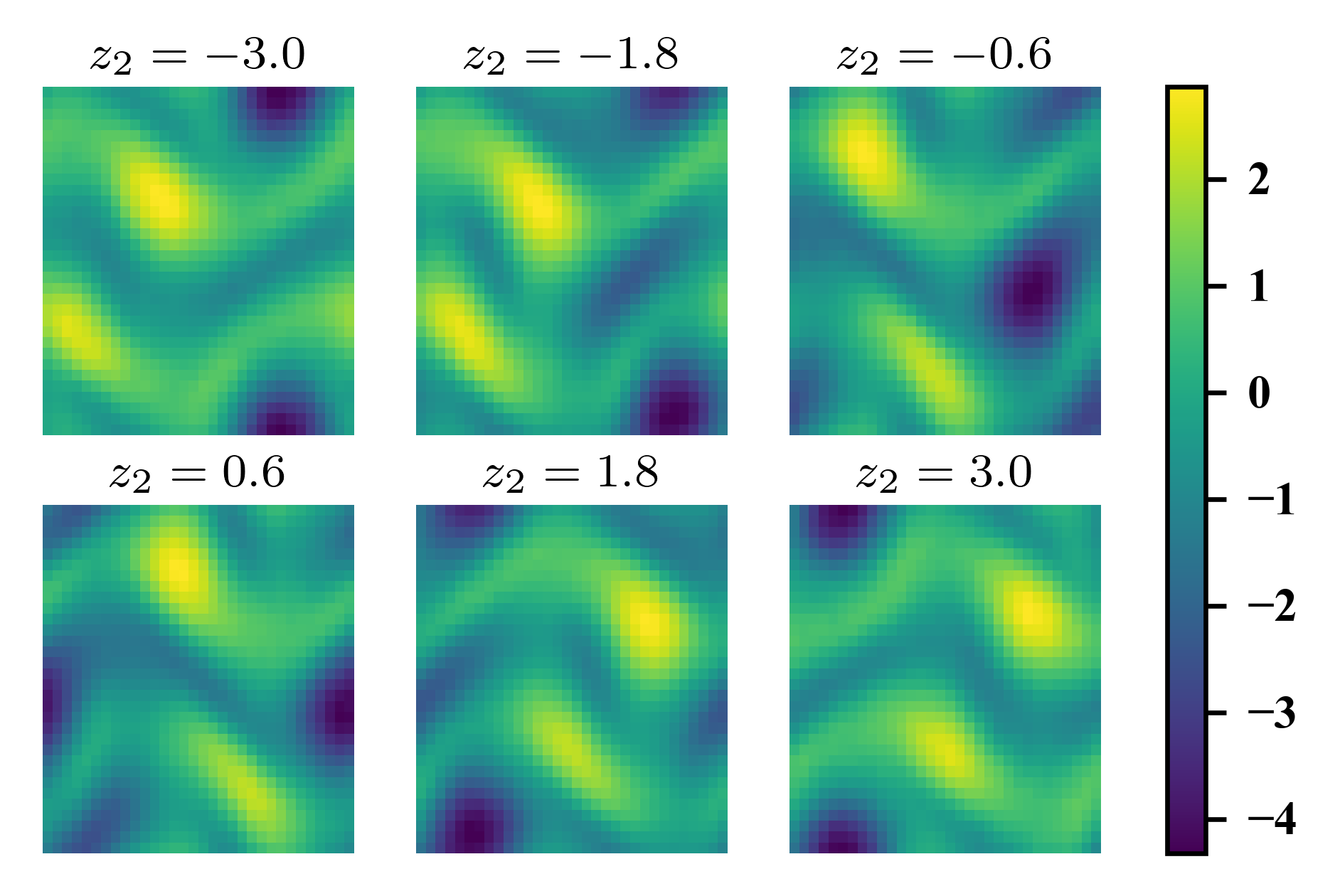}
    \label{fig:z2}
    }
    } \\
    \mbox{
    \subfigure[Reconstructed snapshots with perturbations in $z_3$]{
    \includegraphics[width=0.4\textwidth]{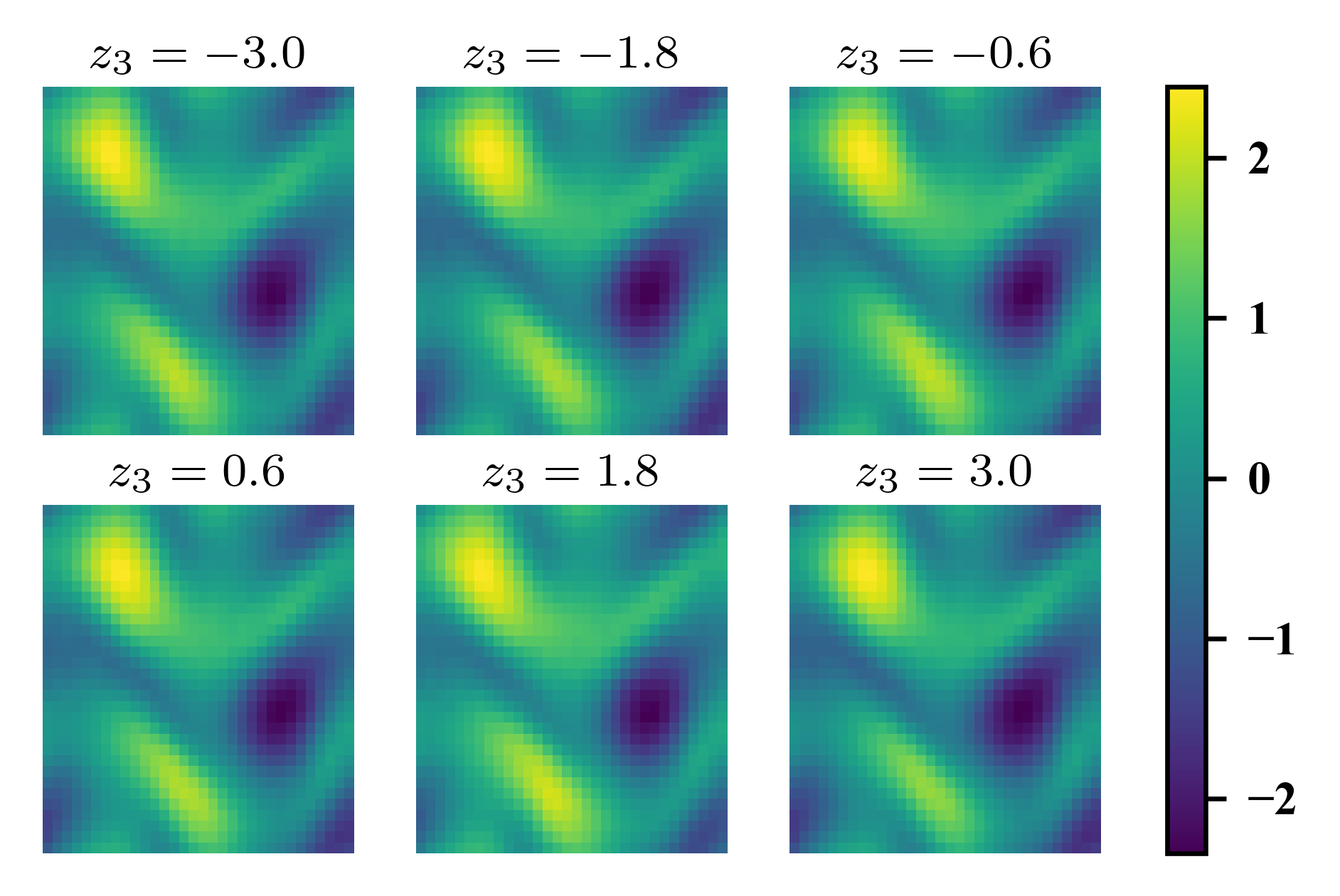}
    \label{fig:z3}
    }
    \subfigure[Reconstructed snapshots with perturbations in $z_4$]{
    \includegraphics[width=0.4\textwidth]{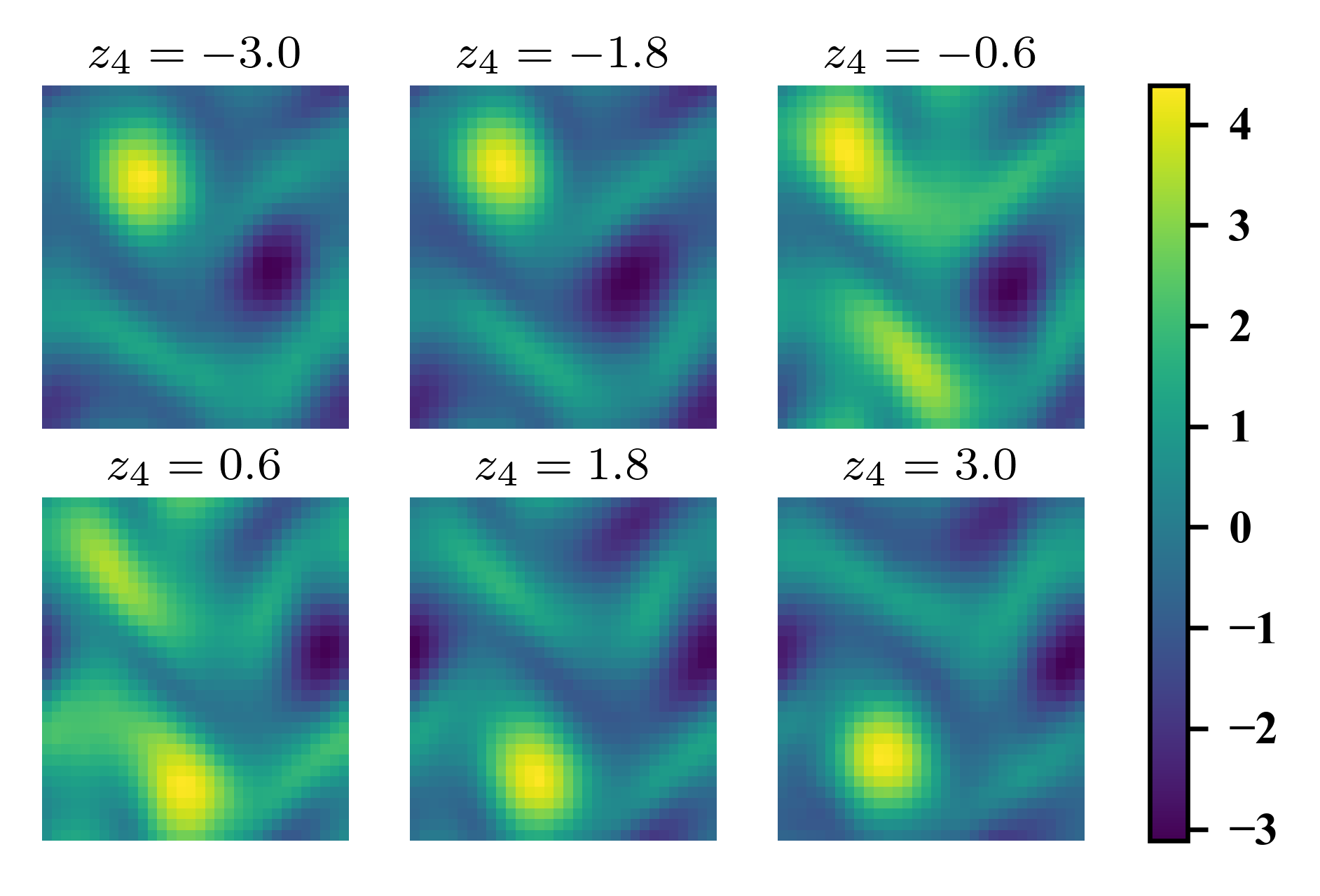}
    \label{fig:z4}
    }
    }
    \caption{Reconstructed snapshots of vorticity from Kolmogorov flow generated by perturbations to latent dimension $z_1, z_2, z_3, z_4$. Perturbations to different states are seen to engender different behavior.}
    \label{fig:disentanglement}
\end{figure*}

We display the reconstruction image of a randomly selected snapshot in Fig.~\ref{fig:original} and Fig.~\ref{fig:reconstructed}, as well as the reconstruction evolution when one dimension of $\boldsymbol{z}$ changes in Fig.~\ref{fig:z1}, Fig.~\ref{fig:z2}, Fig.~\ref{fig:z3}, Fig.~\ref{fig:z4}. The former demonstrates how the latent space information has compressed the high-dimensional state without sacrificing information in the reconstructed snapshot. The latter demonstrates the disentanglement feature of the $\beta-$TCVAE where changing one dimension of $\boldsymbol{z}$ leads to different modes of change in the reconstructed snapshots. For instance, $z_3$ perturbations are seen to have minimal effect on the reconstructed state, while $z_1$ and $z_2$ control translation in the reconstructed snapshot. In constrast, $z_4$ controls reflections about the horizontal plane as the perturbation is moved between positive and negative values.

\subsection{Bayesian online changepoint detection}

In this section, we outline the results from deploying the BOCD algorithm using the latent state obtained by the $\beta$-TCVAE. Given temporal dependencies between the disentangled latent-space time series, we adopt a method to reduce auto-correlation by computing the difference between consecutive time steps, that is, $\boldsymbol{z}_t-\boldsymbol{z}_{t-1}$, which is equivalent to a first-order discretization of the derivative $dz/dt$. This ensures that the assumption of a Markovian process is appropriate. This quantity is utilized for the purpose of changepoint detection. We call this new quantity $\Delta \boldsymbol{z}_t$ (for each dimension of the multidimensional latent state) and furthermore, it is centralized according to 
\begin{align}
    \Delta \boldsymbol{z}_t = \frac{T\Delta \boldsymbol{z}_t -\sum^T_{t=1}\Delta \boldsymbol{z}_t}{\sum^T_{t=1}(\Delta \boldsymbol{z}_t - \frac{1}{T}\sum^T_{t=1}\Delta \boldsymbol{z}_t)^2}
\end{align}

Our latent space, due to the disentanglement property ensures that correlations between dimensions are close to $0$. However, for the purpose of completeness, we still deploy the changepoint detection algorithm while modeling the correlation (i.e., the joint covariance matrix) explicitly for the purpose of computational comparisons. Our results show that the extra modeling costs burdens our computation without bringing extra benefits. We compare the computing time for the BOCPD algorithm for both the independence assumption between latent state dimensions as well through an explicit covariance matrix update as shown in Table \ref{computing}.  In our table, Model 1 refers to the deployment where latent dimensions are assumed to be uncorrelated as a result of disentanglement. Each dimension is modeled with a univariate distribution. In contrast, in Model 2, the latent dimensions are modeled with a multivariate distribution with an online update of the covariance matrix. It is obvious that the use of Model 2 dramatically increases the cost of the changepoint detection algorithm. We note that both models observe a historical period of truncated length $500$.

\begin{table}[!htb]
\centering
\caption{Computing times (in seconds) for Bayesian online change point detection algorithm at different timesteps $T$ of a four-dimensional time series.}
\label{computing}
\begin{tabular}{|c|r|r|}
\hline
$T$ & Model 1 & Model 2 \\ \hline
1000 & 1.02 & 73.86 \\ \hline
2000 & 2.76 & 288.53 \\ \hline
3000 & 5.89 & 656.30 \\ \hline
4000 & 7.52 & 1218.33 \\ \hline
\end{tabular}
\end{table}

\begin{figure*}[!ht]
    \centering
    \mbox{
        \subfigure[Change point detection result with $\alpha=0.01,\ \beta=0.1$, and $\lambda=10$.]{
        \includegraphics[width=0.49\textwidth]{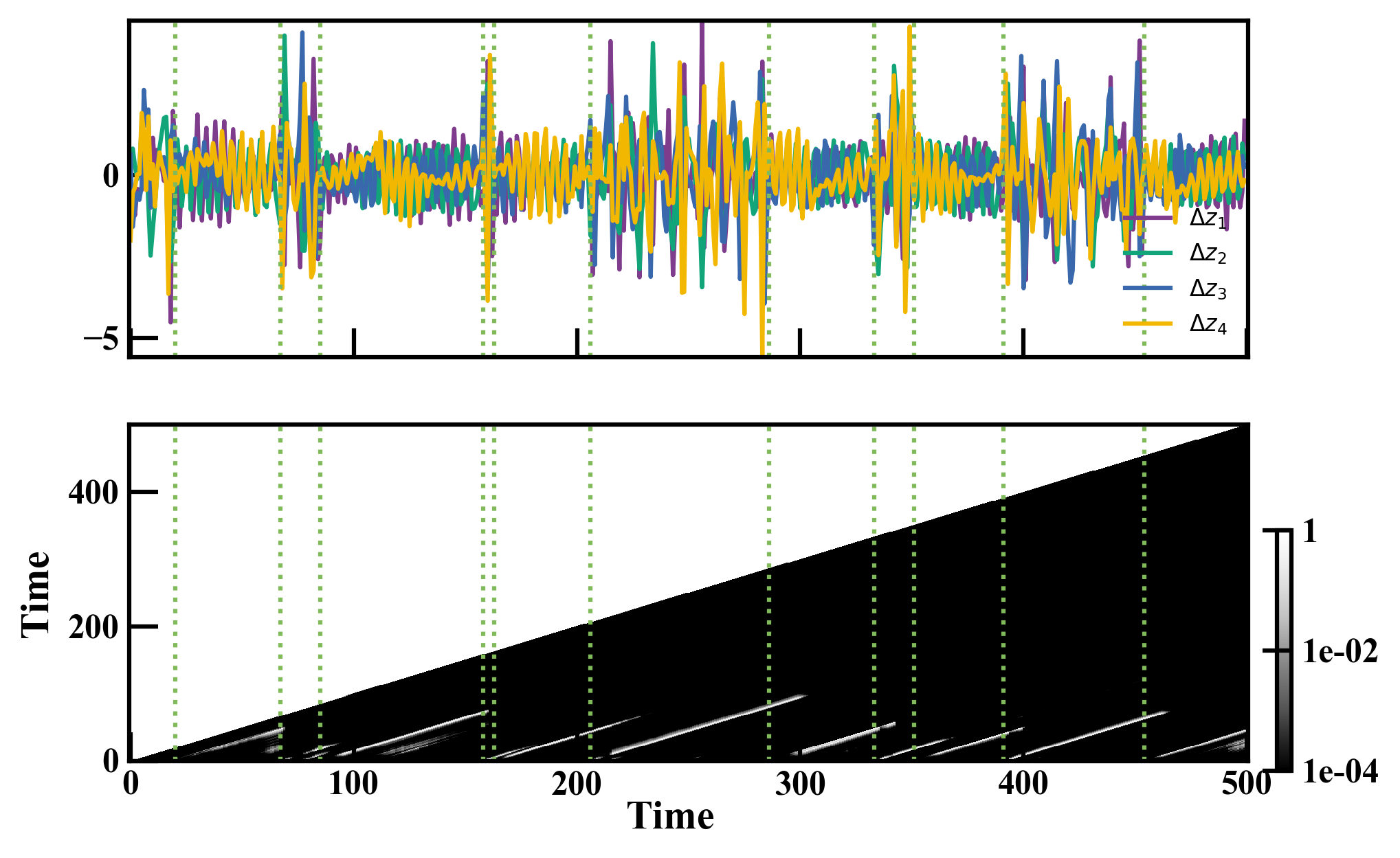}
        \label{cpd1}
        }
        \subfigure[Change point detection result with $\alpha=0.01,\ \beta=0.1$, and $\lambda=50$.]{
        \includegraphics[width=0.49\textwidth]{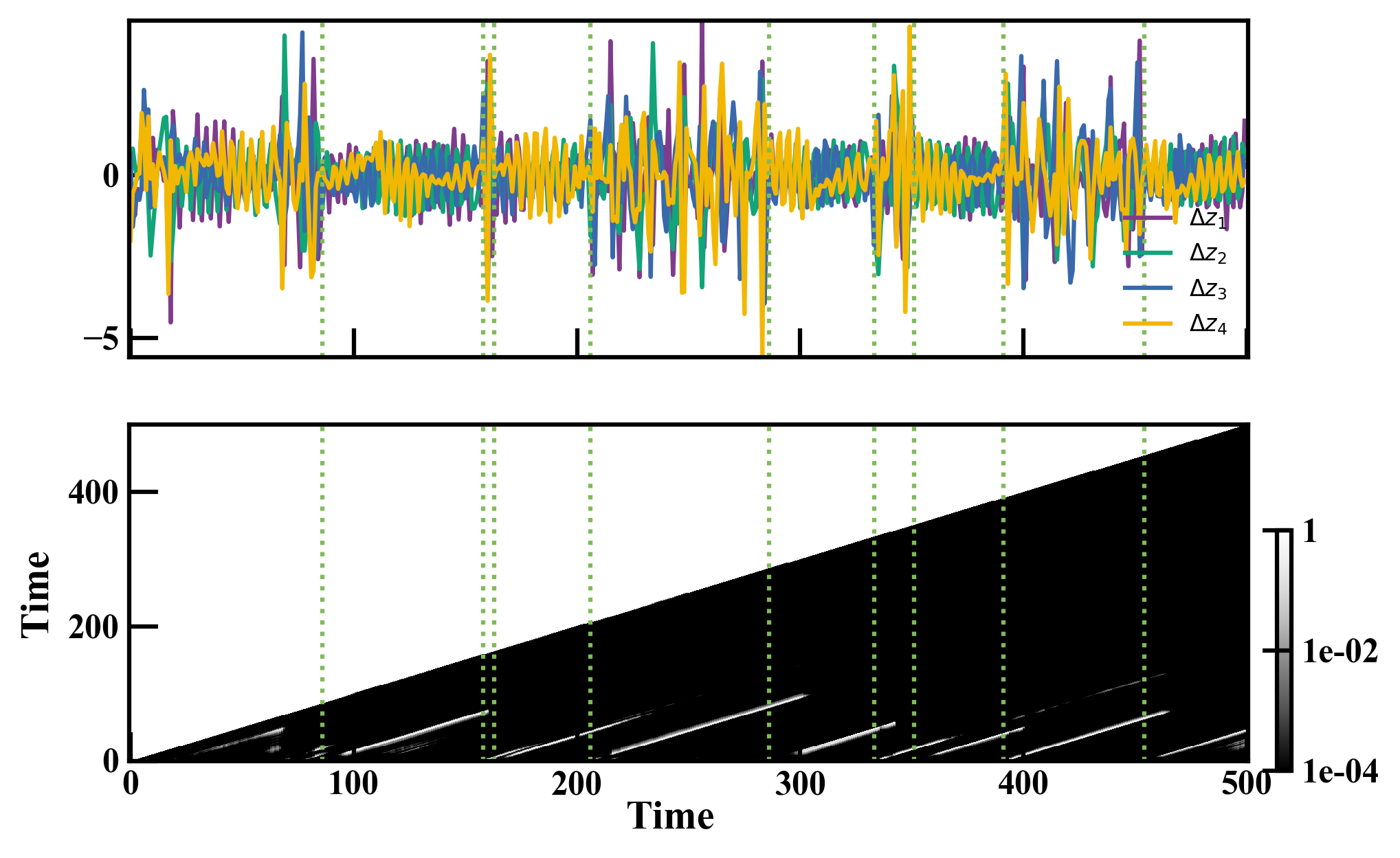}
        \label{cpd2}
        }
    } 
    \caption{The first graph displays the 4-dimensional latent time series $\Delta z$; the right bottom triangle in the second graph represents the probability of a run length from $t_i$ to $t_j$, the horizontal axis represents $t_i$ and the vertical axis represents $t_j$, shallow color represents large probability. Both graphs share the same horizontal axis given by time indices. The green dashed lines are changepoints detected by BOCPD.}
\end{figure*}

As an unsupervised learning algorithm, BOCPD consistently yields reasonable detection outcomes. However, its performance is found to be somewhat sensitive to hyperparameters, as evidenced in Fig. \ref{cpd1} and Fig. \ref{cpd2}, where we show the results of varying the value of $\lambda$. In the first detection task, we assume an expected time interval of $\lambda = 10$ between changepoints, while in the second task, we assume an interval of $\lambda = 50$. Although the difference in the detected changepoints between the two settings is insignificant, the probability of run length calculated by the algorithm exhibits slight variations. In practice, this implies that the most probable changepoints set scheme has some robustness to hyperparameters. However, we still recommend a thorough assessment of hyperparameters based on \emph{a priori} validation with historical data.

\begin{figure}[!h]
 \includegraphics[width=0.45\textwidth]{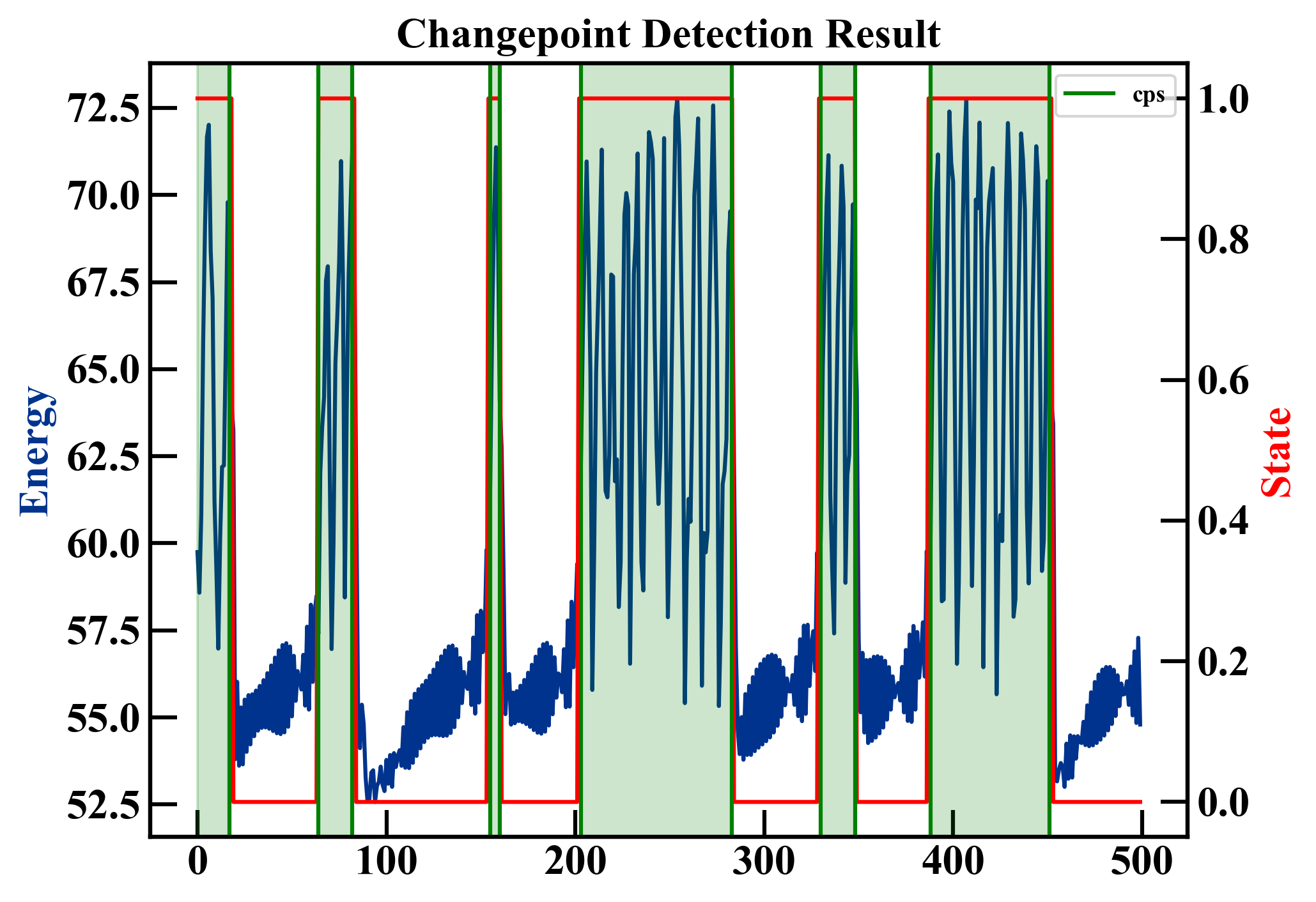}
 \caption{\label{cp}Change point detection results overlapped with kinetic energy. The blue curve is kinetic energy. The red curve is the label of a bursting event provided by human expertise (a value of 1 indicates a bursting event). The green regions represent the difference between two changepoints as detected by our method. The detected change points are provided by the BOCD algorithm with hyperparameter set $C=3$, $\alpha=0.01,\ \beta=0.1$, and $\lambda=10$}
 \label{fig:cpd_result}
\end{figure}

In addition, we found changepoints detected in the latent space $\Delta z$ possesed a delay compared to changepoints provided by manual expertise. Intuitively, this can be assumed to be due to the requirement of observing a critical number of points from the streaming time-series data before the underlying generative process parameterization is updated significantly. For practical applications, this delay can be modeled as an additional hyperparameter during the tuning process with validation data. Here we introduce an additional hyperparameter, $C$, which allows for backward shifting of the detected changepoint to better align with its actual occurrence in the time series, as given by manual expertise. We found $C=3$ in this experiment. In other words, there is a delay of 3 snapshots between the emergence of the bursting event and the detection from the BOCD algorithm looking at its disentangled latent space representation. After tuning such hyperparameters, we achieved optimal detection performance with $\alpha_0=0.01$, $\beta_0=0.1$, and $\lambda=10$. To further ascertain the value of our changepoint detection, we plot our detected points as they overlap with changes in the kinetic energy, which provides the summary statistic of a bursting event as shown in Fig. \ref{fig:cpd_result}. Here. we compare the predicted changepoints with those obtained by manual observations (assumed to be ground truth). We find metrics given by $P=100\%,\ R=100\%,\ F = 100\%$. The comparison indicates that our approach provides accurate detection of the changepoints.

\subsection{Supervised learning based changepoint detection}

Following our BOCD studies, we investigae detecting changepoints as anomalous predictions from a pre-trained time-series forecasting model, given by an LSTM neural network. Our LSTM prediction model was trained using a short time-series of length $2500$, immediately preceding the portion of the time series used for anomaly detection. The LSTM prediction model is tasked with making a prediction for the latent representation of the new snapshot before it is observed, for the purpose of comparison. Before training this model, we prepare a dataset for time series prediction by utilizing the standard sliding window technique. This involves dividing the time series data into overlapping windows, where each window contains a fixed number of time steps. Specifically, we choose a window size of $k$, which determines the number of time steps in each window, and then slide the window over the entire time series data, creating overlapping input-output pairs at each step. For example, if we have a time series with $2500$ time steps and a window size of $10$, we would create $2489$ input-output pairs. Each input-output pair consists of a window of input snapshots and a corresponding window target output snapshots. In practice, this implies a delay of $k$ before a transition detection can be detected in bursting event. Therefore, we prioritize smaller values of $k$. After creating the input-output pairs, we split the dataset into training and validation sets and utilize it to train the LSTM. We train the model for $2000$ epochs with the $L1$ norm as the loss function and apply the model on an unseen test set, where model predictability is used to gauge the emergence of a bursting event. 

\begin{figure*}[!ht]
    \centering
    \mbox{
        \centering
        \subfigure[$k=5$]{
        \includegraphics[width=0.32\textwidth]{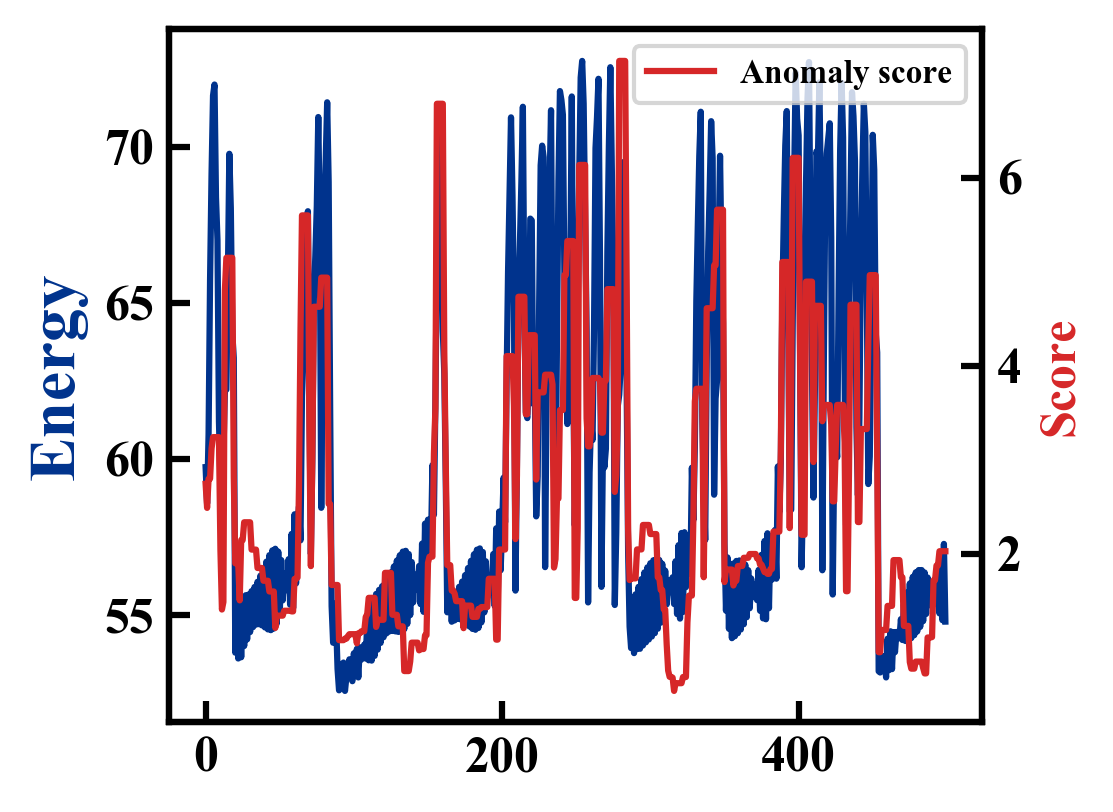}
        \label{fig:anom1}
        }
        \subfigure[$k=10$]{
        \centering
        \includegraphics[width=0.32\textwidth]{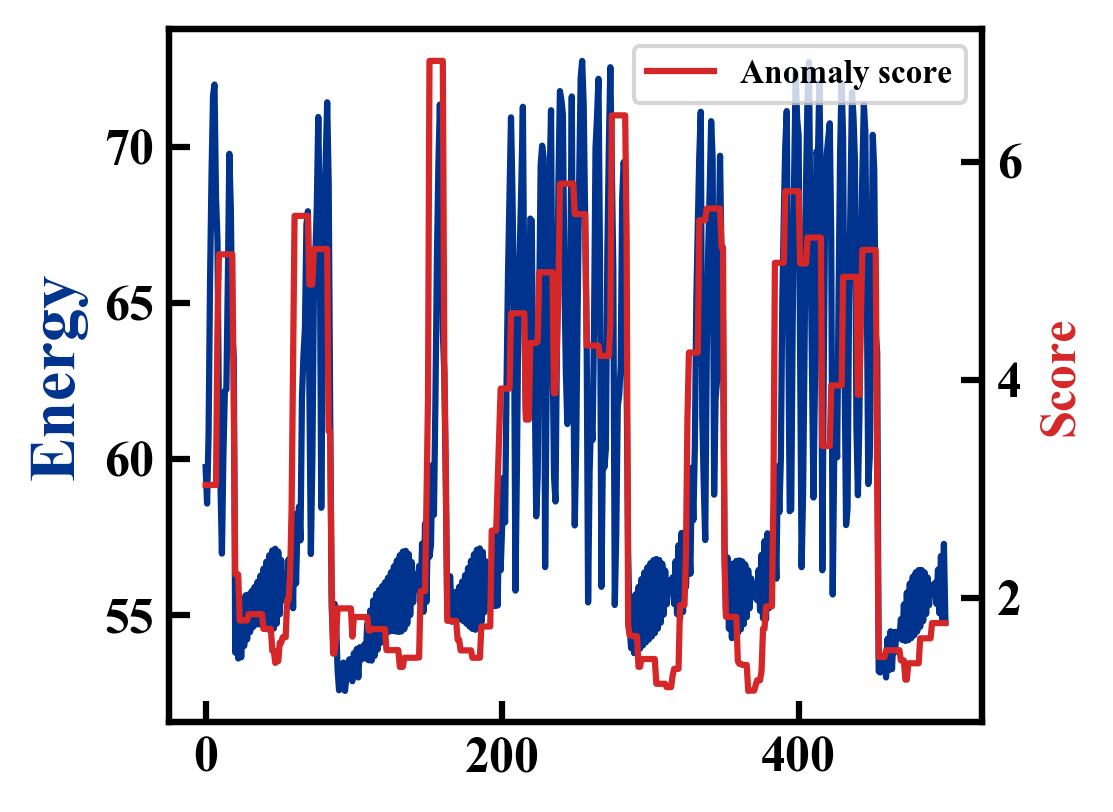}
        \label{fig:anom2}
        }
        \subfigure[$k=15$]{
        \includegraphics[width=0.32\textwidth]{anomaly_score2.png}
        \label{fig:anom3}
        }
    }
    \caption{Anomaly score curves for different LSTM models of varying window lengths.}
\end{figure*}

We conduct a series of experiments to evaluate the performance of LSTM models with varying window lengths, $k$, for detecting anomalies in the Kolmogorov flow dataset. Specifically, we trained three LSTM models and computed anomaly scores for their predictions (see Equation \ref{anom_score}). Here we note that that our observed anomaly score curve (FIG.\ref{fig:anom1}, FIG.\ref{fig:anom3}, FIG.\ref{fig:anom2}) is correlated with the curve that represents the kinetic energy of the Kolmogorov flow and bursting events are detected by the same.

\begin{figure*}[!ht]
    \centering
    \mbox{
        \subfigure[$ k=5, thr=2$.]{
        \centering
        \includegraphics[width=0.33\textwidth]{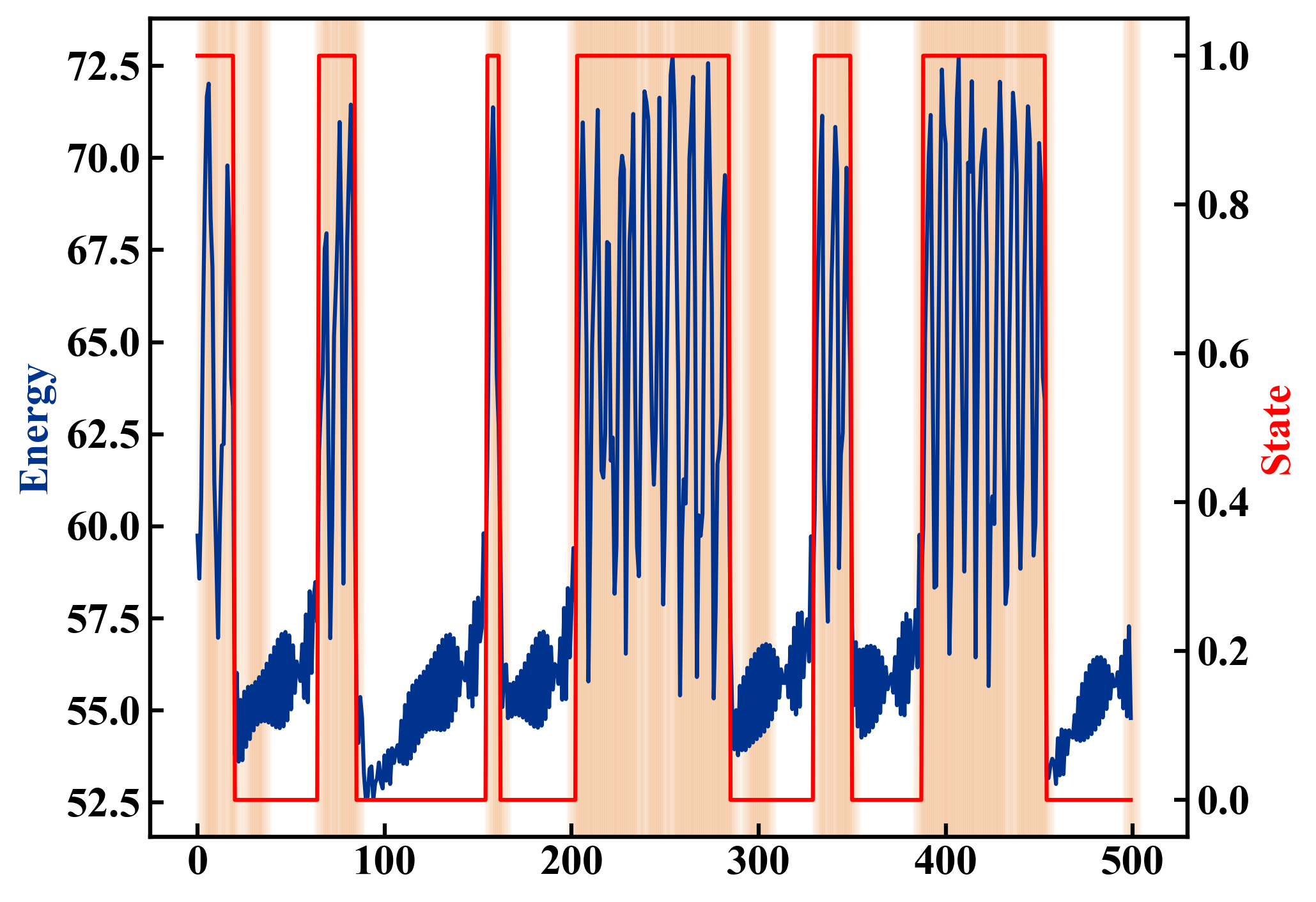}
        \label{fig:anomd11}
        }
        \subfigure[$k=5, thr=2.6$]{
        \centering
        \includegraphics[width=0.33\textwidth]{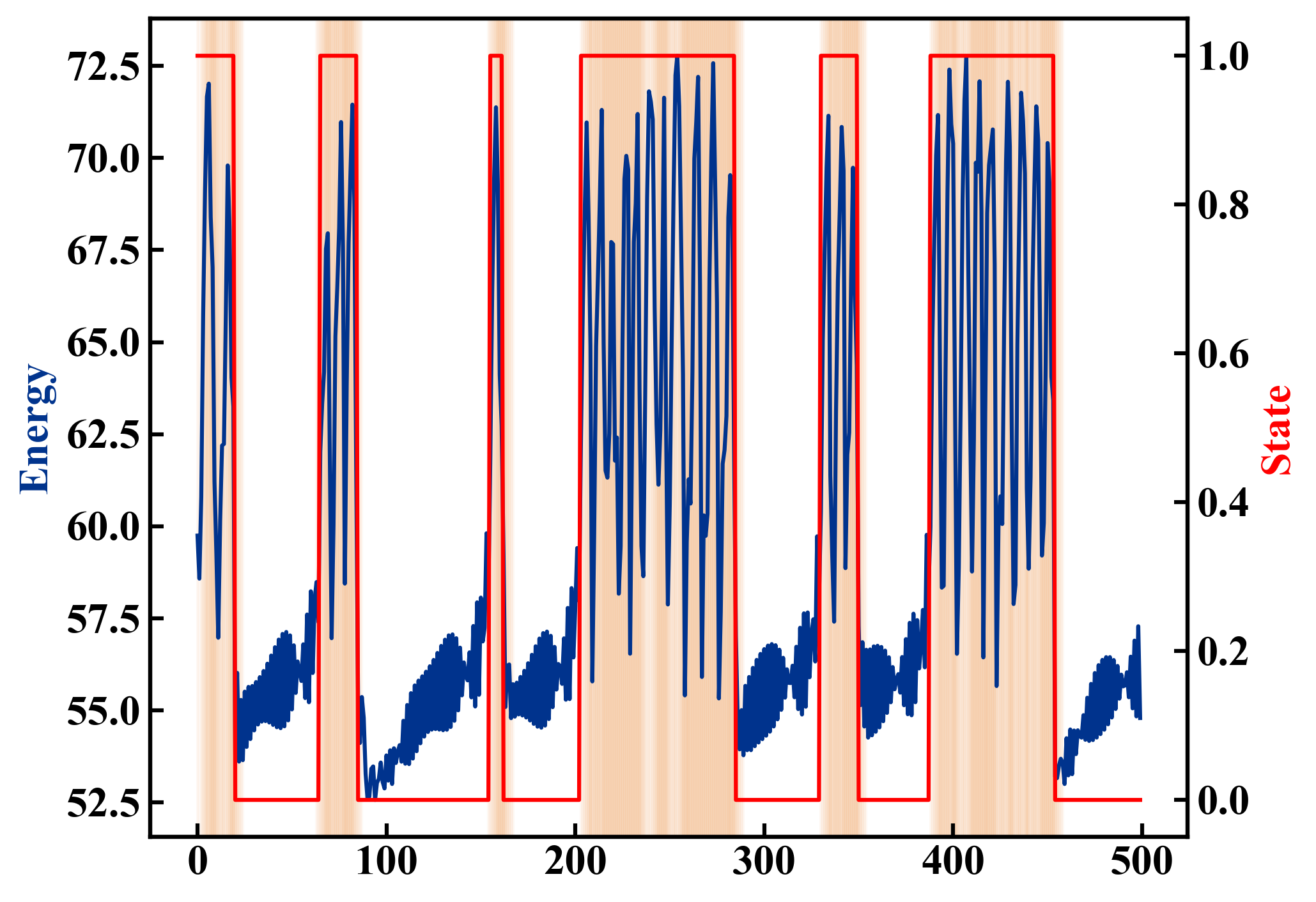}
        \label{fig:anomd12}
        }
        \subfigure[$k=10, thr=2$]{
        \includegraphics[width=0.33\textwidth]{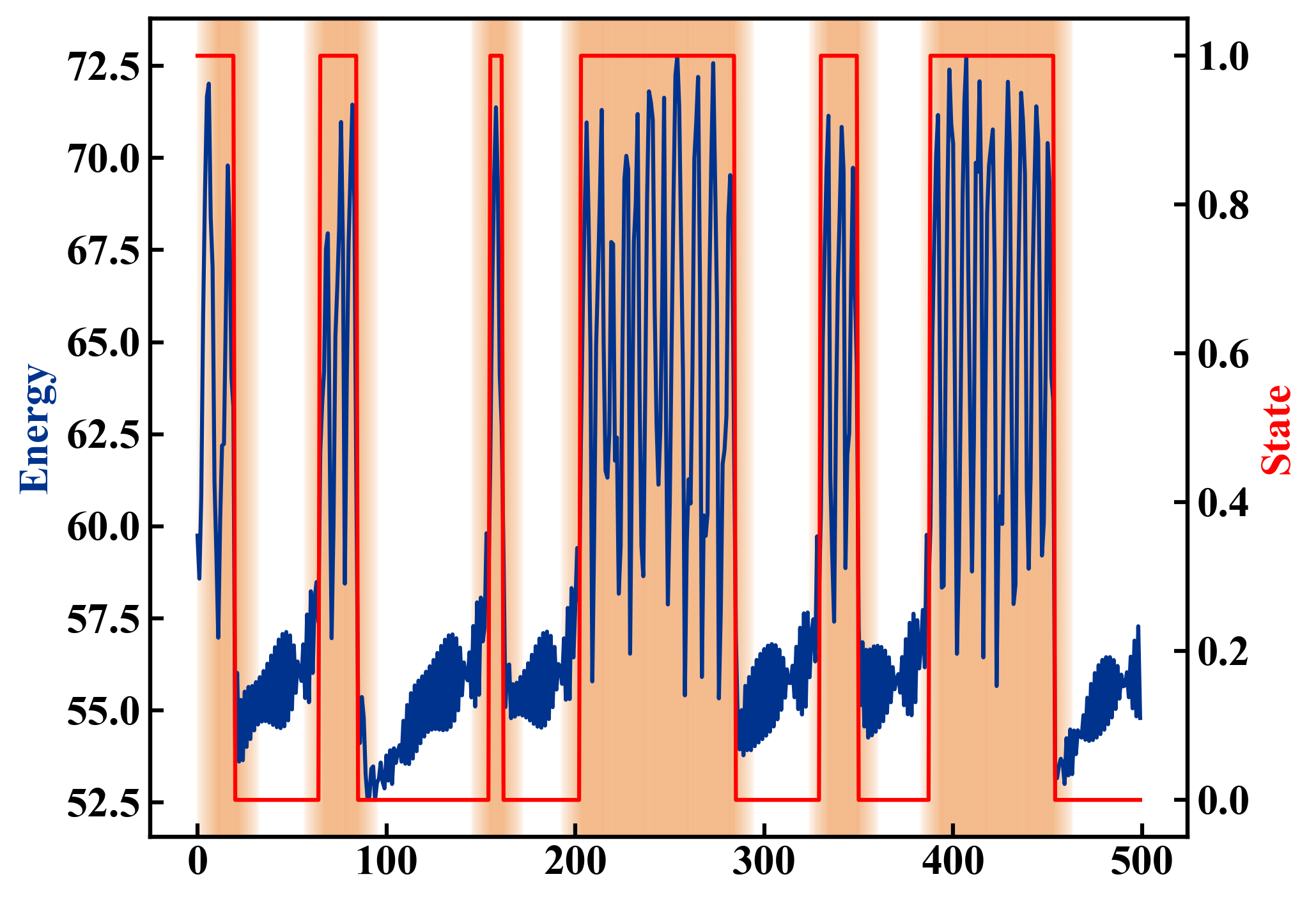}
        \label{fig:anomd21}
        } 
        } \\
    \mbox{
        \subfigure[$k=10, thr=2.6$]{
        \includegraphics[width=0.33\textwidth]{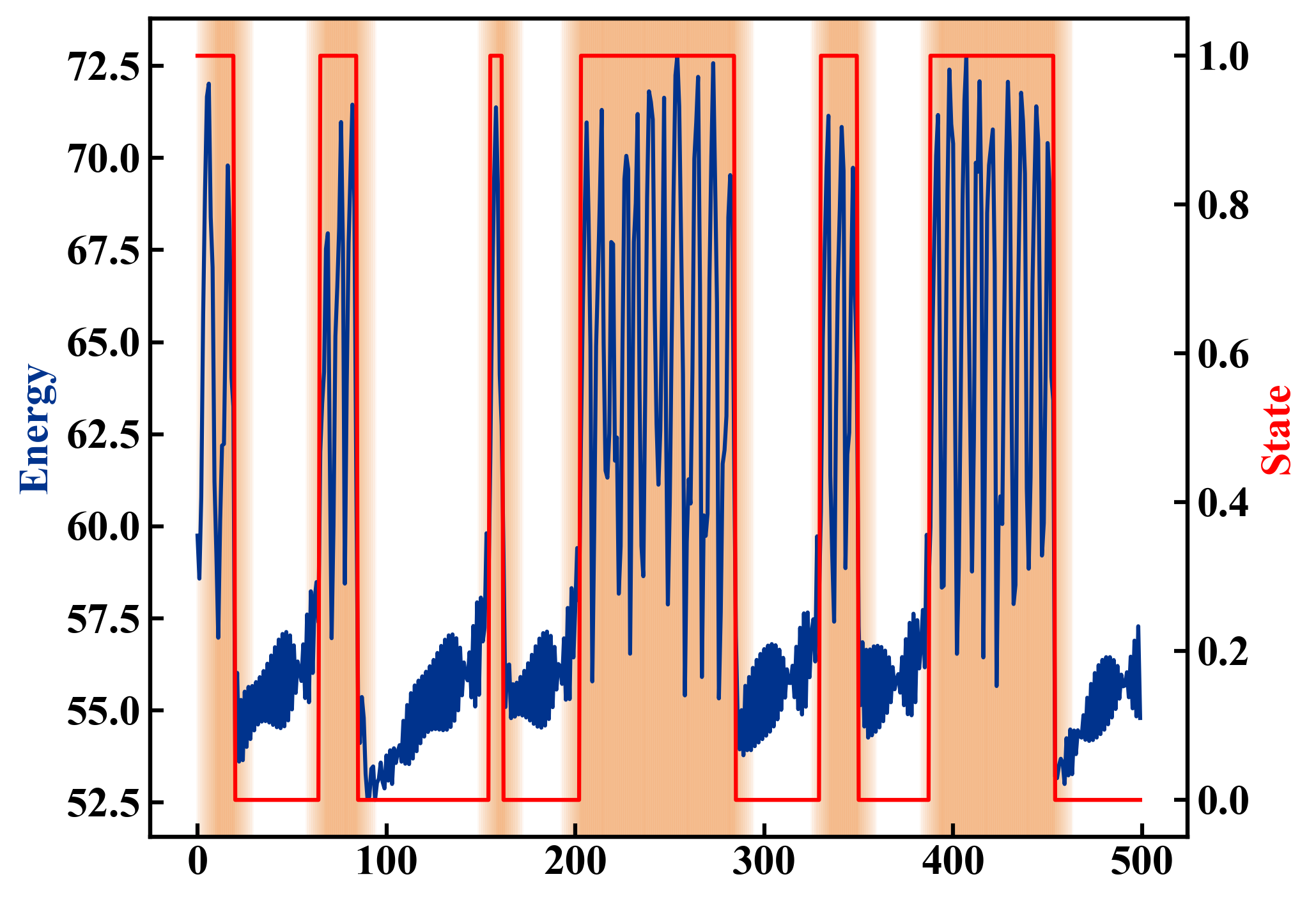}
        \label{fig:anomd22}
        }
        \subfigure[$k=15, thr=2$]{
        \includegraphics[width=0.33\textwidth]{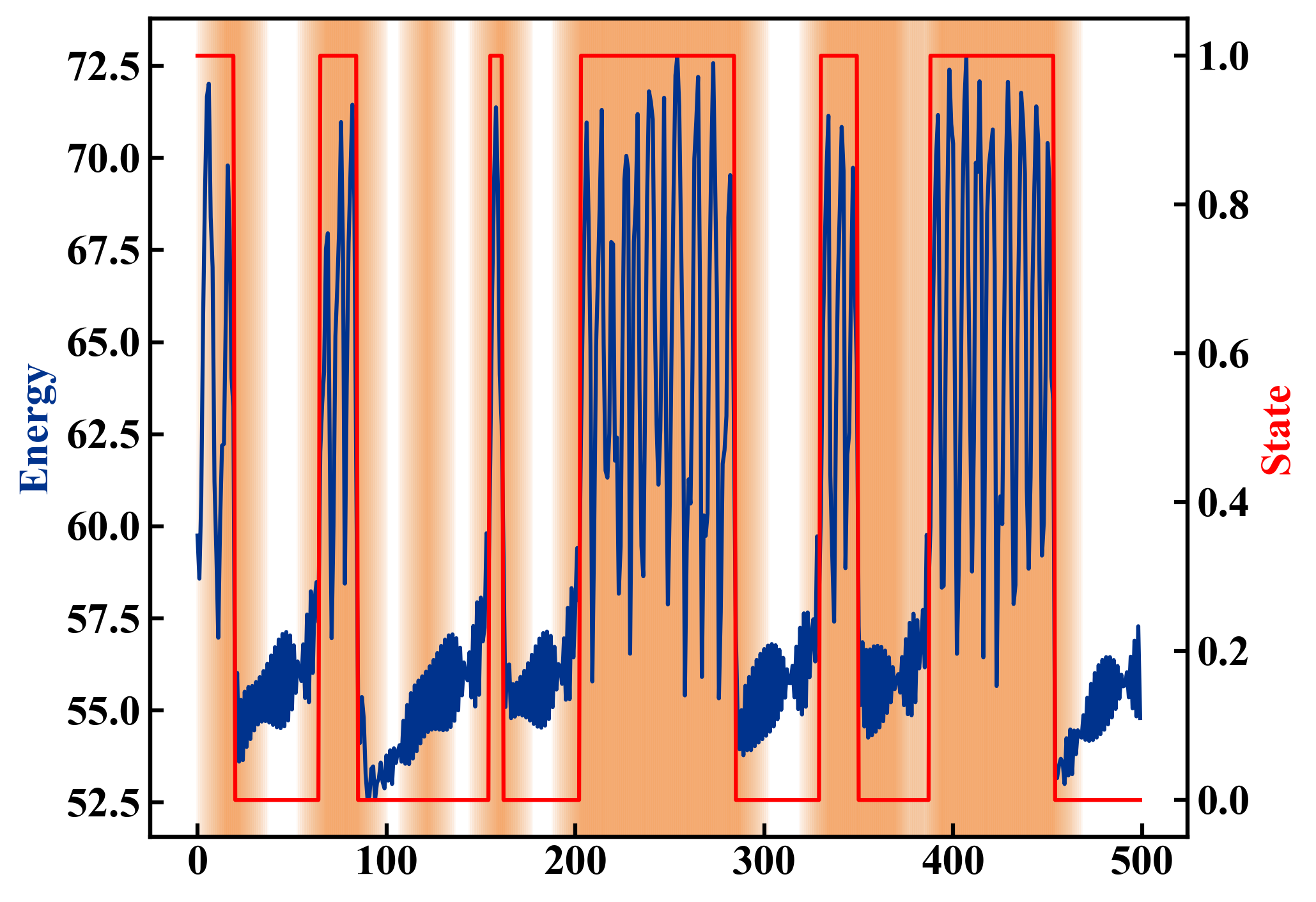}
        \label{fig:anomd31}
        }
        \subfigure[$k=15, thr=2.6$]{
        \includegraphics[width=0.33\textwidth]{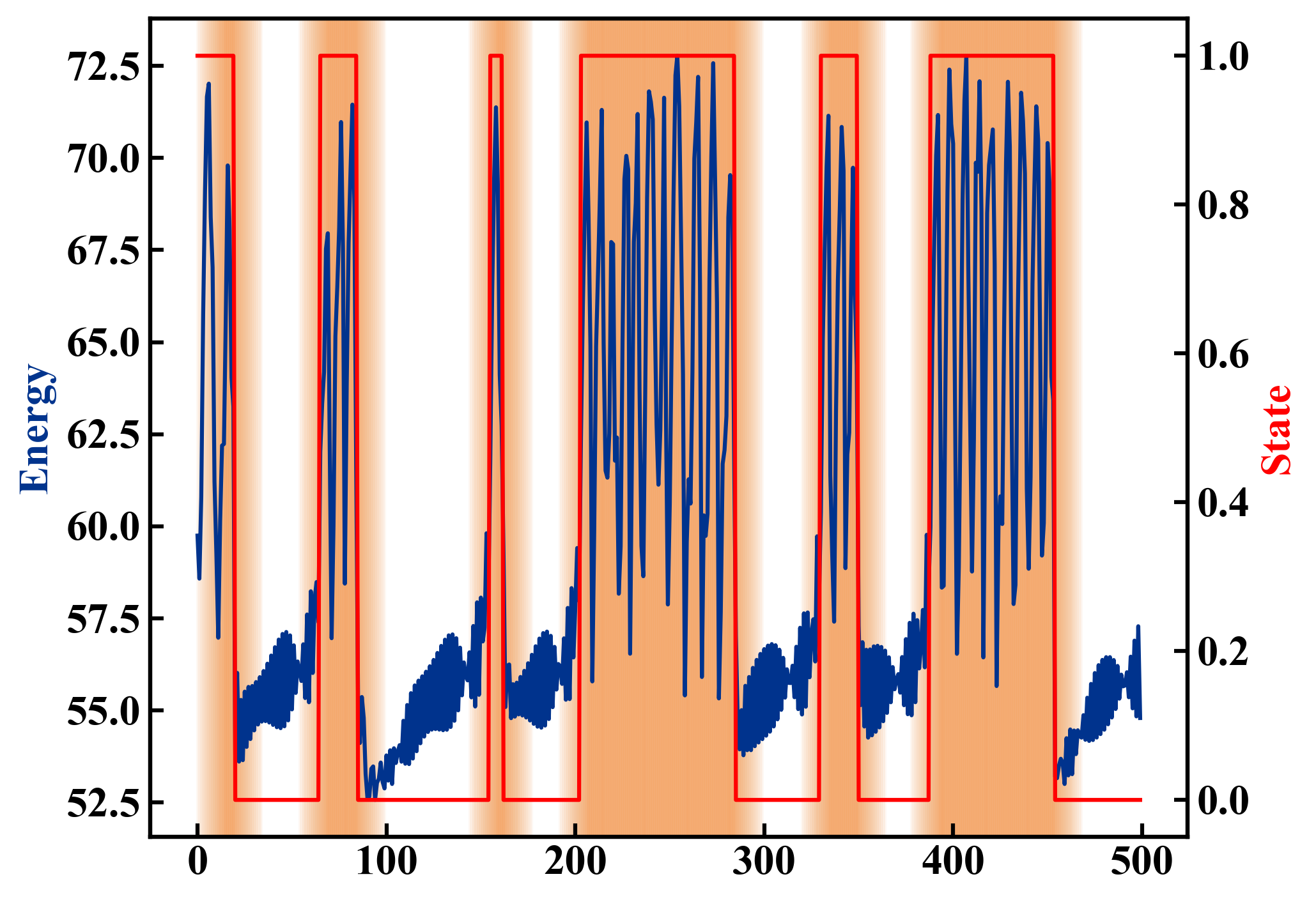}
        \label{fig:anomd32}
        }
    }
    \caption{Changepoint detection results from different LSTM models of varying window lengths, $k$, and anomaly score thresholds $thr$.}
    \label{fig:anomd}
\end{figure*}

Subsequently we explore each model's changepoint detection performance by adjusting the threshold and comparing with the $F-$score. The $F-$score here is the commonly used metric to evaluate classification performance without considering the error margin. Through our experiments, we are able to identify appropriate models and threshold value to achieve accurate anomaly detection using the latent space representation of the Kolmogorov flow dataset. The performance of our different models with varying $k$ and anomaly score threshold $thr$ are shown in Figure \ref{fig:anomd}. A good correlation with the kinetic energy trace is also observed.

Finally, after introducing different techniques to perform changepoint detection using a compressed representation of the original snapshot of vorticity from Kolmogorov flow, we perform an intercomparison using the classical $F-$score. Here, we manually flag classification between changepoints detected by BOCD to calculate the classical $F-$score. In addition, we also identify changepoints based on an optimal threshold of the kinetic energy. This identification technique simply identifies the optimal magnitude of the kinetic energy that closely identifies bursting events.

for the purpose of comparisons. The results are shown in Table \ref{tab:intercomp}.
\begin{table}[!h]
\centering
\footnotesize
\begin{tabular}{|c|c|c|c|c|c|}
\hline
Model & LSTM $k=5$  & LSTM $k=10$ & LSTM $k=15$ & BOCD   & Simple threshold \\ \hline
$F-$Score & 92.78\% & 93.88\% & 89.21\% & 98.33\% & 94.56\% \\ \hline
\end{tabular}
\caption{Anomaly detection result comparison between different methods.}
\label{tab:intercomp}
\end{table}

\section{Conclusion}\label{sec:conclusions}

In this paper, we introduce two complementary techniques for detecting changepoints in high-dimensional dynamical systems. The former is given by a probabilistic method which dynamical adapts the probability of a changepoint given historical information in a Bayesian formulation. The latter is given by a supervised learning framework that flags a changepoint based on inaccurate predictions given training data, i.e., if the truth is significantly different from the prediction, a changepoint is flagged. To accelerate the computation of changepoints, we use a beta variational autoencoder to identify a latent space that reconstructs the high-dimensional snapshots of our system, with the added benefit of independence between coordinates in the latent space. By enforcing independence, we observe significant speed-up in the deployment of the probabilistic technique mentioned above. We also note that the probabilistic approach, while interpretable and deployable in an online fashion, requires careful selection of hyperparameters through a potentially costly search. Therefore, given a section of exemplar data for a dynamical system including its changepoints, the probabilistic and deep learning based approaches are both competitive for detecting regime changes. Finally, we remark that both techniques explored in this article may be enhanced by using state-of-the-art techniques such as universal density approximators (normalizing flows and diffusion models) instead of a conjugate-exponential models and transformer-based approaches instead of the LSTM. We also recognize challenges for the proposed methods to generalize to situations where the changepoint detection is desired \emph{a priori}. Our follow-up studies will explore such considerations. 

\section*{Acknowledgements}

This work was supported by the U.S. Department of Energy (DOE), Office of Science, Office of Advanced Scientific Computing Research (ASCR), under Contract No. DEAC02–06CH11357, at Argonne National Laboratory. We also acknowledge funding support from ASCR for DOE-FOA-2493 ``Data-intensive scientific machine learning''. SL acknowledges the Givens-Associates program at Argonne National Laboratory for enabling this work.

\bibliography{maipsamp}

\end{document}